\documentclass{amsart}

\usepackage{amsfonts, amsmath, amsthm}
\usepackage{epsfig, psfrag}
\newtheorem{theorem}{Theorem}[section]

\theoremstyle{definition}

\theoremstyle{remark}
\newtheorem{remark}[theorem]{Remark}

\numberwithin{equation}{section}

\swapnumbers \theoremstyle{plain}

\newtheorem{thm}{Theorem}[section]
\newtheorem*{thm2}{Theorem}
\newtheorem{lem}[thm]{Lemma}

\newtheorem{cor}[thm]{Corollary}

\newtheorem{prop}[thm]{Proposition}

\theoremstyle{remark}
\newtheorem*{rem}{Remark}

\theoremstyle{definition}

\newcommand{\T}{\mathcal T}
\newcommand{\C}{\mathcal C}

\renewcommand{\P}{{\mathcal P}}
\renewcommand{\S}{{\mathcal S}}

\newcommand{\bdy}{\partial}
\newcommand{\bbb}{\mathbb}
\newcommand{\rppp}{\mathbb{R}P^3}
\newcommand{\rpp}{\mathbb{R}P^2}

\newcommand{\td}{\tilde}
\newcommand{\open}[1]{\stackrel{\circ}{#1}}

\newcommand{\abs}[1]{\lvert#1\rvert}

\begin{document}
\title{Finding planar surfaces in knot- and link-manifolds }
\author{William Jaco}
\address{Department of Mathematics, Oklahoma State University,
Stillwater, OK 74078}

\email{jaco@math.okstate.edu}
\thanks{The first author was partially supported by NSF Grants
 DMS9971719, DMS0204707 and DMS0505609,
The Grayce B. Kerr Foundation, The American Institute of Mathematics
(AIM), and The Lois and Fred Gehring Visitor Chair at University of
Michigan}

\author{J.~Hyam Rubinstein}
\address{Department of Mathematics and Statistics,
University of Melbourne, Parkville, VIC 3052, Australia}
\email{rubin@maths.unimelb.edu.au}
\thanks{The second author was partially supported by The Australian Research
Council and The Grayce B. Kerr Foundation.}

\author{Eric Sedgwick}
\address{School of CTI,
DePaul University, 243 S. Wabash Ave, Chicago, IL 60604}
\email{esedgwick@cti.depaul.edu}
\thanks{The third author was partially supported by US-Israel BSF2002039 and The Grayce B. Kerr Foundation.}


\date{\today}


\begin{abstract} It is shown that given any link-manifold, there is an
algorithm to decide if the manifold contains an embedded, essential
planar surface; if it does, the algorithm will construct one. The
method uses normal surface theory but does not follow the classical
approach. Here the proof uses a re-writing method for normal
surfaces in a fixed triangulation and may not find the desired
solution among the fundamental surfaces. Two major results are
obtained under certain boundary conditions. Given a link-manifold
$M$, a component $B$ of $\bdy M$, and a slope $\gamma$ on $B$, it is
shown that there is an algorithm to decide if there is an embedded
punctured-disk in $M$ with boundary $\gamma$ and punctures in $\bdy
M\setminus B$; if there is one, the algorithm will construct one.
Again, while normal surfaces are used, we may not find a solution
among the fundamental surfaces. In this case we use induction on the
number of boundary components of the link-manifold. It also is shown
that given a link-manifold $M$, a component $B$ of $\bdy M$, and a
meridian slope $\mu$ on $B$, there is an algorithm to decide if
there is an embedded punctured-disk with boundary a longitude on $B$
and punctures in $\bdy M\setminus B$; if there is one, the algorithm
will construct one. This is shown to follow from the previous result
using a link-manifold related to $M$ and called the link-manifold
obtained from $M$ by Dehn drilling along the slope $\mu$. The
properties of minimal vertex triangulations, layered-triangulations,
$0$--efficient triangulations and especially triangulated Dehn
fillings are central to our methods. We also use an average length
estimate for boundary curves of embedded normal surfaces; the
average length estimate shows, in quite general situations, that
given the link-manifold $M$ by a triangulation $\T$, then all normal
surfaces of a bounded genus must have a short boundary curve on some
boundary of $M$. The constant that determines how short is
completely determined by the fundamental surfaces in $(M,\T)$. A
version of the average length estimate with boundary conditions also
is derived.

\end{abstract}

\maketitle

\section{introduction}

This work began in the eighties with an attempt to develop a
singular normal surface theory as a means toward solving the Word
Problem for $3$--manifold groups. The Word Problem for
$3$--manifolds can be formulated as a decision problem given a
knot-manifold.

A compact, orientable $3$--manifold with nonempty boundary, each
component of which is a torus, is called a {\it link-manifold}. If
the boundary is connected, then we say it is a {\it knot-manifold}.
We are interested in algorithms to determine if a given knot- or
link-manifold contains an interesting planar surface. Often we are
interested in how the planar surface sits within the manifold and
for this we use some special terminology.

The isotopy class of a non contractible simple closed curve in a
torus is called a {\it slope}. Classically, the study of
knot-manifolds has been as the exterior of knots embedded in some
other manifold. In such a situation, there is a unique slope on the
boundary of the knot-manifold, corresponding to the isotopy class of
curves on its boundary that bound a disk in the solid torus
neighborhood of the knot. Such a curve is called a meridian. When we
are given a knot-manifold along with a slope in a component of its
boundary, we use the term {\it meridian} for the given slope. A
slope with geometric intersection {\it one} with the meridian is
called a {\it longitude}. There are infinitely many longitudes for a
meridian, each obtained from any of the others by Dehn twists about
the meridian.

If $D$ is a disk, $\{p_1,\ldots,p_K\}$ are distinct points in the
interior of $D$,  and $\{\eta(p_1),\dots,\\\eta(p_K)\}$ are pairwise
disjoint regular neighborhoods of the points with $\eta(p_i)\subset
\open{D}, 1\le i\le K$, then we say that
$P=D\setminus\bigcup_{i=1}^K\open{\eta}(p_i)$ is a {\it
punctured-disk with boundary $bd(P)=\bdy D$ and punctures
$\bdy\eta(p_1),\ldots,\bdy\eta(p_K)$}. Of course, a punctured-disk
is a planar surface but a punctured-disk has a distinguished
boundary component, and all other boundary components are called
punctures.

In this setting, the following is the Word Problem for the
fundamental groups of closed $3$--manifolds.

\vspace{.125 in}\noindent {\footnotesize WORD PROBLEM} (closed
$3$--manifolds). {\it Given a knot-manifold $M$ and a meridian on
$\bdy M$. Decide if a longitude bounds a (possibly) singular
punctured-disk in $M$ with punctures meridians}.

\vspace{.125 in}We note that if any longitude bounds a singular
punctured-disk, then all longitudes bound a singular punctured-disk.

Our approach was to understand singular punctured-disks by
considering them as normal surfaces. One quickly notices that the
analogous question for an embedded punctured-disk, which must
actually be a disk having no punctures, is the Classical Unknotting
Problem. The Unknotting Problem was solved for knot-manifolds in
$S^3$ by W. Haken \cite{haken1}, where he showed that given a knot
in $S^3$ it can be decided if it is the unknot. On the other hand,
there is a curious analogy to link-manifolds, which does not have so
fortunate an outcome; in fact, it is an insoluble problem.

\vspace{.125 in}\noindent {\footnotesize WORD PROBLEM} (finitely
presented groups). {\it Given a link-manifold $M$, a component $B$
of $\bdy M$, and a meridian slope on $B$. Decide if there is a
(possibly) singular punctured-disk in $M$ with boundary slope a
longitude in $B$ and punctures in $\bdy M\setminus B$ or meridians
on $B$}.

\vspace{.125 in} The Word Problem for finitely presented groups in
not solvable; hence, the preceding is an insoluble decision problem
for link-manifolds.

In the last section of this paper, we discuss the equivalence of
these statements of the  Word Problem with other familiar versions.

Being unable to make progress in the case of singular normal
surfaces, we decided to investigate if analogous problems for
embedded surfaces had solutions and found several interesting
questions regarding embedded planar surfaces in knot- and
link-manifolds. Unfortunately, at that time, we did not make much
progress on these either and laid the problems aside until the late
nineties. In the late nineties, we discovered a number of new tools
for working with normal surfaces and special triangulations, both of
which lend themselves nicely to algorithmic problems. We soon
obtained that given a link-manifold it can be decided if there is an
embedded, essential planar surface in the manifold. We use the term
{\it essential} in describing a properly embedded surface in a
$3$-manifold in this work to mean that the surface is incompressible
and is not parallel into the boundary. Recall that for
link-manifolds, an incompressible surface is also
$\bdy$--incompressible or is an annulus; and if the knot-manifold is
irreducible, the annulus must be parallel into the boundary. It was
at that time we returned to these problems. We delayed writing these
results for publication until now.

In Section 4 we give the following general result about finding
interesting planar surfaces.

\begin{thm2} Given a link-manifold there is an algorithm to decide if
it contains a properly embedded, essential, planar surface; if it
does, the algorithm will construct one.\end{thm2}

Possibly the most interesting aspects of this result are the tools
used and the method of its proof. Typically, algorithms using normal
surface theory follow a standard format.

Firstly, there is an existence step where one shows that if a given
manifold contains an embedded surface with a property $\P$, then it
contains a normal one with this property $\P$. However, at some
points in this work, we need to  modify the triangulation in order
to assure the existence of a normal surface with the desirable
property.

Secondly, there is a recognition step. The necessary algorithms for
recognition that a given normal surface is an essential surface are
given in \cite{jac-tol}; however, we also use from \cite{jac-tol}
that if $F$ is an  embedded, essential surface and is least weight
in its isotopy class, then every normal surface that has its
projective representative in the carrier of $F$ also is an embedded,
essential, normal surface. The latter result first appeared in
\cite{jac-oer} using handle-decompositions and later in
\cite{bart-sch} for closed surfaces and triangulations.

Thirdly, in the classical approach, the big (typically, by far the
hardest) step is to show that if there is a surface with a property
$\P$, then there is one among the fundamental surfaces.

We note that a positive solution to the first two steps places us in
the common situation for recursively enumerable problems. Namely,
given a $3$--manifold $M$ by a triangulation $\T$, all normal
surfaces in $M$ (with respect to $\T$) can be constructed. Thus we
merrily go about constructing the normal surfaces. If one with
property $\P$ exists, it is normal and if it is normal, we can
recognize it. Thus if the given manifold has such a surface, we will
eventually find one. However, if there is none, we do not know this
and do not know when to stop looking. There are only finitely many
fundamental surfaces; hence, if our surface must be among this
finite set, which we can construct straight away, then we have
solved the problem. Our algorithms do not have the classical step
three; in fact, we might need to go quite far afield of the
fundamental solutions of the given triangulation but we do find a
finite set in which to look and all surfaces in this finite set are
normal in our given triangulation. In Section 4, the proof is by a
method we call a {\it re-writing process}. A re-writing process was
also used in \cite {jac-rub0}; the re-writing process here is
different but is the same principle.

In Section 5, we answer the problem for embedded punctured-disks
analogous to the singular problem above for the Word Problem for
finitely presented groups. It is given in the second of the next two
theorems.

\begin{thm2} Given a
link-manifold $M$, a component $B$ of $\bdy M$, and  a slope
$\gamma$ in $B$, there is an algorithm to decide if $M$ contains an
embedded punctured-disk with boundary having slope $\gamma$ and
punctures in $\bdy M\setminus B$. If there is one, the algorithm
will construct one.\end{thm2}

\begin{thm2} Given a link-manifold $M$, a component $B$ of
$\bdy M$, and a meridian in $B$, there is an algorithm to decide if
$M$ contains an embedded punctured-disk with boundary having slope
of a longitude in $B$ and punctures in $\bdy M\setminus B$. If there
is one, the algorithm will construct one.\end{thm2}

Again, we do not conform to the classical third step of finding our
solution among the fundamental surfaces. Moreover, our proof uses an
interesting method for link-manifolds; we use induction on the
number of boundary components. For rather subtle reasons, we can not
use induction in Section 4.

Besides not using the classical form of proof for the above
problems, we also call upon a number of new results on
triangulations and new tools in normal surface theory. We discuss
what we need from the literature and provide results necessary for
this work in Section 3 and later sections. For example, our methods
typically require minimal-vertex triangulations or at least
triangulations that have at most one vertex in each boundary
component. In some situations we need $0$--efficient triangulations,
which are minimal vertex triangulations for knot- and
link-manifolds. These triangulations are quite general and given a
$3$--manifold via any triangulation, there are algorithms that
modify the given triangulation to one of these that fits into our
methods. In particular, in Section 4, we generalize a result from
\cite{jac-rub0} and prove the following prime decomposition theorem,
where $n(Q)$ in the connected sum decomposition means the connected
sum of $n$ copies of the manifold $Q$ and $Card(\T)$ for a
triangulation $\T$ stands for the number of tetrahedra of $\T$.

\begin{thm2} Given a link-manifold $M$ via a triangulation
$\T$, there is an algorithm to construct a prime decomposition
$$M=p(S^2\times S^1)\# q(\rppp)\# r(D^2\times S^1)\#
M_1\#\cdots\#M_n,$$ where $p,q$ and $r$ are nonnegative integers and
each $M_i$ is given by a $0$--efficient triangulation $\T_i,  i
=1,\ldots,n$, respectively; furthermore, $\sum_{i=1}^n Card(\T_i)\le
Card(\T)$.\end{thm2}

We also use results from \cite{jac-rub-layered} and
\cite{jac-sedg-dehn} on layered-triangulations of the solid torus
and the classification of normal surfaces in minimal layered
triangulations of the solid torus. This is used in conjunction with
triangulated Dehn fillings, which were introduced in these same two
references. If we are given a link-manifold $M$ with a triangulation
$\T$, having just one vertex in each boundary component, then for
slopes $\alpha_1,\ldots\alpha_K$ on distinct boundary components of
$M$, there is a natural way to triangulate the Dehn filling
$M(\alpha_1,\ldots\alpha_K)$ with a triangulation that is $\T$ on
$M$ and is a minimal layered-triangulation of each of the solid tori
added. We discuss this and some of the results we use from the
literature in Section 3.

We assume the reader is familiar with the basic concepts from normal
surface theory; however, in Section 3, we identify some of the
particular concepts and results that we use. Among these, one which
may not be so familiar, is that mentioned above where if $F$ is an
embedded, essential normal surface and is least weight in its
isotopy class, then every normal surfaces that projects into the
carrier of $F$ is embedded and essential. We also define and use the
length of the boundary of a normal surface. Another major tool
related to the boundary of a normal surface, also used in
\cite{jac-sedg-dehn}, is what we call the A{\it verage} L{\it ength}
E{\it stimate} (ALE) for the boundary of a normal surface. For a
manifold $M$ with triangulation $\T$, there is a constant $C$
dependent only on $M$ and $\T$ so that every normal surface of
bounded genus has the average length of its boundary no large than
$C$; hence, in particular, for link-manifolds (under the right
conditions) a normal planar surfaces must have a short boundary
slope on some boundary of the link-manifold. This enables us to find
a finite family of (short) slopes in which to do Dehn fillings and
apply an induction hypothesis. In using triangulated Dehn fillings
in conjunction with ALE, we are able to have very strong control on
the complexity of our methods.

Among the questions we consider is one to determine if there is a
planar normal surface with a prescribed boundary slope on some
boundary component. This and other questions that involve the
boundary of the normal surface are called {\it boundary conditions}.
Section 5 considers a theory of normal surfaces with boundary
conditions. In particular, we show that for the triangulations we
use, given a slope on the boundary, there is a set of matching
equations so that we get a normal solution space having only normal
surfaces that meet the given boundary component in the given slope.
We may have closed normal surfaces and surfaces that meet other
boundary components in a totally uncontrolled way. These method can
be greatly generalized; we did not do that here. We also adapt ALE
to normal surfaces with boundary conditions in Section 5.

In the last section we return to a discussion of those versions of
the Word Problem given above and what we might call the classical
versions. Following the solution of the Geometrization Conjecture,
we know that the Word Problem for $3$--manifold groups is solvable;
however, we remain curious as to the existence of a straight forward
method, say, in the spirit of the solution to the Word Problem for
the fundamental groups of Haken manifolds, given by F. Waldhausen
\cite{wald-word}.

\section{Background Material } We shall assume the reader has familiarity with our notion of
triangulations, as well as a basic knowledge of normal surface
theory. The references \cite{jac-rub0, jac-sedg-dehn,
jac-rub-layered} serve as good background material for both
triangulations from our point of view and basic facts on normal
surfaces. We have, however, collected some facts and background in
this section, which are particularly relevant to this work.

\subsection{Triangulations.} One of the interesting aspects of this
work is not only the effectiveness of minimal vertex triangulations
in understanding the combinatorial and algorithmic problems we
encounter but also in enabling us to understand the topology better.
The following is typical of the type of triangulations we like.

\begin{thm}\cite{jac-sedg-dehn,jac-rub-blowup}\label{minvertex-triangs} Suppose
$M$ is a compact, orientable $3$--manifold with boundary, no
component of which is a $2$--sphere. Then any triangulation of $M$
can be modified to a triangulation having all vertices in $\bdy M$
and just one vertex in each component of $\bdy M$.\end{thm}

A proof is given in both the cited references. The idea is quite
straight forward. First, there is an algorithm due to R.H. Bing
\cite{bing1} that can be used to modify a given triangulation to one
having all vertices on the boundary. This step was not mentioned but
is necessary in the proof given in \cite{jac-sedg-dehn}. Having all
vertices in the boundary, then the method of ``closing-the-book",
described in detail in Theorem 3.3 of \cite{jac-sedg-dehn},
finishes the proof.

A manifold having a triangulation with just one vertex (closed or
possibly compact and bounded with just one boundary component) or a
triangulation with all vertices in the boundary and just one vertex
in each boundary component (compact with boundary) can not have a
triangulation with fewer vertices. We shall refer to such
triangulations as {\it minimal-vertex triangulations}. This is
different from a minimal triangulation of the manifold $M$; a
triangulation $\T$ of $M$ is a {\it minimal triangulation} if and
only if for any triangulation $\T'$ of $M$, $Card(\T)\le Card(\T')$,
where $Card(\T)$ is used to denote the number of tetrahedra of a
triangulation. For irreducible $3$-manifolds, distinct from $S^3,
\rppp$, $L(3,1)$, and $\bbb{B}^3$, minimal triangulations are
minimal-vertex triangulations. We suspect minimal triangulations are
also minimal-vertex triangulations for reducible $3$--manifolds but
have not established this. Minimal-vertex triangulations are
completely general; the above theorem shows that any compact
$3$--manifold with boundary (no component of which is a $2$--sphere)
admits a minimal-vertex triangulation and a proof is given in
\cite{jac-rub-layered} that any closed $3$--manifold admits a
one-vertex triangulation.

\subsection{Normal surfaces.} We collect here the main results from normal surface
theory we will need along with our conventions for notation. Primary
sources for this material are
\cite{jac-ree,jac-rub4,jac-sedg-dehn,jac-rub0,jac-tol}.

If $M$ is a $3$--manifold, a triangulation of $M$ selects a family
of surfaces called normal surfaces. Typically a normal surface is
defined as an embedded surface that meets the tetrahedra of the
triangulation in normal triangles and normal quadrilaterals. An
isotopy of $M$ that is invariant on the various simplices of the
triangulation is called a {\it normal isotopy}. A normal isotopy
class of normal triangles or of normal quadrilaterals is called a
{\it triangle type} or {\it quad type}, respectively. We caution the
reader that with the triangulations we are using, where the
simplicies are only embedded on their interiors and may have
identifications on their boundaries, an embedded surface is normal
if and only if its pull back to the tetrahedra before face
identifications is a collection of normal triangles and normal
quadrilaterals.

We shall use standard normal coordinates. In this case, there are
four normal triangle types and three normal quad types for each
tetrahedron, giving each normal surface a parametrization with $7t$
variables, where $t$ is the number of tetrahedra in the
triangulation. Hence, an embedded normal surface determines a unique
nonnegative integer lattice point in $\bbb{R}^{7t}$. The
triangulation also determines a system of homogeneous linear
equations, the {\it matching equations}; its  solution space meets
the nonnegative orthant of $\bbb{R}^{7t}$ in a cone called the {\it
solution cone}, which we denote $\S(M,\T)$.  We say two normal
surfaces satisfy the same {\it quadrilateral condition} if they do
not meet any one tetrahedron in distinct quadrilateral types;
equivalently, there is an additional set of conditions placed on the
coordinates in the solution cone where for each tetrahedron two of
the quadrilateral types have been set to zero. This algebraic
condition is also called a quadrilateral condition. There are $3^t$
possible quadrilateral conditions. Each integer lattice point in the
solution cone corresponds to a (possibly singular) normal surface.
This correspondence is one-one between embedded normal surfaces in
$M$ with respect to $\T$ and the integer solutions in the solution
cone that also satisfy a quadrilateral condition. Those points in
the solution cone that have norm one in the $\ell_1$--norm ($\sum
x_i = 1, x_i\ge 0$) form a compact, convex, linear cell called the
{\it projective solution space} for $M$ with triangulation $\T$,
which is denoted $\P(M,\T)$. The solution cone is the cone over
$\P(M,\T)$ with vertex the origin. Each point in the solution cone
has a unique projection into projective solution space. If two
points in the solution cone project to the same point in $\P(M,\T)$,
we say they are {\it projectively equivalent}. If $F$ is a point in
the solution space and $\overline{F}$ is its projection in
$\P(M,\T)$, then we call the minimal dimensional closed face of
$\P(M,\T)$ containing $\overline{F}$ the {\it carrier of $F$} and
denote it by $\C(F)$. We will not distinguish notation between an
embedded normal surface $F$ and its parametrization $F$ in
$\bbb{R}^{7t}$.

From the Hilbert Basis Theorem, there is a {\it unique, minimal}
finite set of integer lattice solutions in $\S(M,\T)$,
$F_1,\ldots,F_K$, so that for $F$ any integer lattice point in
$\S(M,\T)$, we have $$F = \sum n_iF_i,$$ where $n_i$ is a
nonnegative integer. We call such a family $F_1,\ldots,F_K$ {\it
fundamental solutions}. If $V$ is an integer lattice point in
$\S(M,\T)$, $V$ projects to a vertex $\overline{V}$ of $\P(M,\T)$,
and if for any integer lattice point $V'$ that projects to
$\overline{V}$ we have $V' = kV$ for some positive integer $k$, we
call $V$ a {\it vertex solution}. All vertex solutions must be among
any set of fundamental solutions. The vertex solutions may be found
using any one of a number of methods from linear programming and the
fundamental solutions may be found once one has the vertex
solutions. We have the following characterizations of fundamental
and vertex solutions.

\begin{enumerate}\item{\it $F$ is a fundamental solution if and only
if whenever $A$ and $B$ are solutions and $F = A + B$,  either $A=0$
or $B=0$.}\item{\it $V$ is a vertex solution if and only if whenever
$A$ and $B$ are solutions and there is a positive integer $k$ so
that $kV = A + B$, then $A = k'V$ and $B = k''V$, $k', k''$ positive
integers.}\end{enumerate}

If $F$ and $F'$ are normal surfaces and they satisfy the same
quadrilateral conditions, then using standard cut-and-paste
techniques from $3$--manifold topology, there is a unique way to
form a normal surface from $F$ and $F'$ called the {\it geometric
sum} of $F$ and $F'$. Since the geometric sum of two embedded normal
surfaces $F$ and $F'$ is parameterized by the coordinate sum of
their parameterizations, we also write the geometric sum of $F$ and
$F'$ as $F+F'$.

There are several forms of complexity associated with normal
surfaces; two of the simplest ones are its weight, analogous to
area, and the length of its boundary. If $F$ is a normal surface in
$M$, then $F$ is in general position with respect to the
$2$--skeleton of the triangulation $\T$; we define
$wt(F)=Card(F\cap\T^{(1)})$ to be the {\it weight of $F$},  where
$Card(S)$ is the cardinality of $S$. Similarly, we define $L(\bdy
F)=Card(\bdy F\cap\T^{(1)})$ as the {\it length of $\bdy F$}.

There is a special form for geometric addition; namely, the
geometric sum $F+G$ is said to be in {\it reduced-form} if for all
possible ways to write $F+G$ the number of components of $F\cap G$
is minimal; i.e., if $F+G = F'+G'$, then the number of components of
$F\cap G$ is no larger than the number of components of $F'\cap G'$.
We have the following very useful observations, which we learned
from \cite{sch}. Clearly, any geometric sum may be written in
reduced-form.

\begin{lem} Suppose the geometric sum $F+G$ is defined and in reduced-form.
Then\begin{itemize}\item[-] a component of $F\cap G$ does not
separate both $F$ and $G$;\item[-] if $F+G$ is connected, then both
$F$ and $G$ are connected.\end{itemize}\end{lem}

We note the following easily established facts where $F$ and $G$ are
embedded normal surfaces that satisfy the same quadrilateral
conditions.

\begin{enumerate}\item $\chi(F + G) = \chi(F) +\chi(G)$,\vspace{.05 in}\item
$wt(F+G) = wt(F) + wt(G)$, {\it and}\vspace{.05 in}\item $L(\bdy
(F+G)) = L(\bdy F) + L(\bdy G)$.\end{enumerate}

\subsection{Other Triangulations.} There are some special
triangulations we will be using and we briefly discuss these
triangulations and related useful facts.

\vspace{.1 in}\noindent {\bf - $\boldsymbol{0}$--efficient
triangulations.} A triangulation of a closed $3$--manifold is said
to be {\it $0$--efficient} if and only if the only normal
$2$--spheres are vertex-linking; if the manifold has boundary, a
triangulation is said to be {\it $0$--efficient} if and only if the
only normal disks are vertex-linking. From \cite{jac-rub0} we have
the following about $0$--efficient triangulations.

Suppose $\T$ is a $0$--efficient triangulation of the $3$--manifold
$M$.

If $M$ is closed, then\begin{enumerate}\item[(i)] $M$ is irreducible
and contains no embedded $\rpp$.  \item[(ii)] $\T$ has one vertex or
$M=S^3$; if $M = S^3$, then $\T$ has at most two
vertices.\end{enumerate}

If $M$ has nonempty boundary, then\begin{enumerate}\item[(iii)] $\T$
has no normal $2$--spheres,\item[(iv)] $M$ is irreducible and
$\bdy$--irreducible. \item[(v)] All the vertices of $\T$ are in
$\bdy M$ and there is just one vertex in each boundary component or
$M$ is a $3$--cell.\end{enumerate}

Hence, a $0$--efficient triangulation is a
minimal-vertex-triangulation, except possibly for $S^3$ and the
$3$--cell. For a $3$--cell with a $0$--efficient triangulation, then
the triangulation is expected to have precisely three vertices, all
in the boundary; it is easy to see all vertices must be in the
boundary but we have not been able to show there are only three.

In Theorem \ref{find-planar}, we show that given a link-manifold
$M$, we can construct a prime decomposition of $M$, where the
irreducible and $\bdy$--irreducible factors have $0$--efficient
triangulations. It is shown in \cite{jac-rub0} that any compact,
irreducible, $\bdy$--irreducible, orientable $3$--manifold, distinct
from $\rppp$, admits a $0$--efficient triangulation. In particular,
minimal triangulations of these manifolds are $0$--efficient.

\vspace{.1 in}\noindent {\bf - One-vertex triangulations and slopes
in tori.} Up to homeomorphism of the torus there is a unique
one-vertex triangulation. It has two triangles, three edges and (of
course) one vertex. For any triangulation of a surface an essential
(not contractible) simple closed curve is isotopic to a normal
curve; however, for a one-vertex triangulation of a torus there is
more. Namely, by Lemma 3.5 \cite{jac-sedg-dehn} essential curves in
a one-vertex triangulation of a torus are isotopic if and only if
they are normally isotopic; thus in such a triangulation there is a
unique normal isotopy class for each essential simple closed curve.
By Lemma 3.4 \cite{jac-sedg-dehn}, we also note that in any
one-vertex triangulation of a closed surface, the only trivial
(contractible) normal curve is vertex-linking.  Hence, in a
one-vertex triangulation of a torus, slopes and normal isotopy
classes of essential simple closed curves are in one-one
correspondence.

In a one-vertex triangulation of a torus, we say two slopes are {\it
complementary} if the geometric sum of their normal representatives
is a union of trivial (vertex-linking) curves. Each slope has a
unique complementary slope. The following (Proposition 3.7 of
\cite{jac-sedg-dehn}) is a fundamental result about slopes of
boundaries of {\it normal} surfaces in $3$--manifolds having tori in
their boundary and triangulations that induce one-vertex
triangulations on these boundary tori.

\begin{thm}\label{slopes}\cite{jac-sedg-dehn} Let $M$ be an orientable $3$--manifold having a component
of its boundary a torus, $T$, and let $\T$ be a triangulation of $M$
that restricts to a one-vertex triangulation of $T$. Suppose $S_1$
and $S_2$ are embedded normal or almost normal surfaces and $\bdy
S_1\subset T$. If $S_1$ and $S_2$ satisfy the same quadrilateral
conditions and both meet $T$ in non-trivial slopes, then these
slopes are either equal or complementary.\end{thm}

Some components of the boundaries of $S_1$ and $S_2$ may be trivial
curves in the boundary of $M$; however, it is implicit in the
theorem that there is an essential curve from each of $S_1$ and
$S_2$ in $T$ to determine slopes. This theorem gives the result that
for a knot-manifold with a triangulation inducing a one-vertex
triangulation on the boundary torus, there are only finitely many
boundary slopes for normal and almost normal surfaces. In
particular, it gives the result, discovered earlier by A. Hatcher
\cite{hatcher-slope}, that for a knot-manifold $M$ there are a
finite number of slopes bounding embedded, incompressible and
$\bdy$--incompressible surfaces in $M$.

If we choose two slopes, having isotopy classes $\lambda$ and $\mu$
on a torus and geometric intersection one, then they determine a
basis for the first homology of the torus. Representing their
homology classes also by $\lambda$ and $\mu$, respectively, we have
that any slope $\alpha$ can be represented as $\alpha =
a\lambda+b\mu$, where $a$ and $b$ are relatively prime integers
(possibly $a=0, b=\pm 1$ or $a=\pm 1,b=0$). We define the {\it
distance} between the slopes $\alpha$ and $\beta$, denoted $\langle
\alpha,\beta\rangle$, to be their geometric intersection number;
hence, for a basis $\lambda, \mu$ and $\alpha = a\lambda + b\mu,
\beta = c\lambda+ d\mu$, we have $\langle \alpha,\beta\rangle =
\abs{ad-bc}$. $\langle \cdot,\cdot\rangle$ is not a true distance
function but we do have $\langle \alpha,\beta\rangle = 0$ if and
only if $\alpha = \beta$ and $\langle\alpha,\beta+\gamma\rangle
=\langle\alpha,\beta\rangle +\langle\alpha,\gamma\rangle$ and
$\langle\alpha,m\beta\rangle=m\langle\alpha,\beta\rangle$. Once we
have designated a basis for the homology of the torus, slopes are in
one-one correspondence with $\bbb{Q}\cup{\infty}, \bbb{Q}$ the
rationals.

\vspace{.1 in}\noindent{\bf - Layered-triangulations of the solid
torus and triangulated Dehn fillings.} Layered-triangulations of the
solid torus are  studied extensively in \cite{jac-rub-layered}. Both
layered triangulations of the solid torus and triangulated Dehn
fillings are used in \cite{jac-rub1,jac-sedg-dehn,jac-rub-layered}.
We provide a brief review here.

Suppose $M$ is a compact $3$--manifold with nonempty boundary, $\T$
is a triangulation of $M$ and $\T_{\bdy}$ is the induced
triangulation on $\bdy M$. Furthermore, suppose $e$ is an edge in
$\T_{\bdy}$ and there are two distinct triangles $\sigma$ and
$\beta$ in $\T_{\bdy}$ meeting along the edge $e$. Let $\td{\Delta}$
be a tetrahedron distinct from the tetrahedra in $\T$ and let
$\td{e}$ be an edge in $\td{\Delta}$. Suppose $\td{\sigma}$ and
$\td{\beta}$ are the faces of $\td{\Delta}$ that meet along
$\td{e}$. We can identify $\td{e}$ with $e$ and extend this to face
identifications from $\td{\sigma}\rightarrow\sigma$ and
$\td{\beta}\rightarrow\beta$, getting a  $3$--manifold $M'$
homeomorphic with $M$ and a triangulation $\T'$ of $M'$, having one
more tetrahedron than $\T$. We write $M' = M\cup_e\Delta$ and $\T' =
\T\cup_e\td{\Delta}$, where $\Delta$ is the image of $\td{\Delta}$
and say that $\T'$ is obtained from $\T$ by  {\it layering (a
tetrahedron) on $\T$ along the edge $e$}. Notice that this operation
transforms  the triangulation $\T_{\bdy}$ to the triangulation
$\T'_{\bdy}$ by what is called a Pachner or bi-stellar move of type
$2\leftrightarrow 2$ on $\T_{\bdy}$ along the edge $e$ (also called
a ``diagonal flip" within the quadrilateral $\sigma\cup\beta$).

\begin{figure}[htbp]

           \psfrag{D}{\Large{$\td{\Delta}$}}
\psfrag{a}{\small{$a$}} \psfrag{b}{\small{$b$}}
\psfrag{c}{\small{$c$}} \psfrag{d}{\small{$d$}}
\psfrag{x}{\small{$x$}}\psfrag{y}{\small{$y$}}\psfrag{z}{\small{$z$}}
\psfrag{M}{\begin{tabular}{c}
           {\footnotesize $\langle
a,b,c\rangle\leftrightarrow\langle x,y,z\rangle$}\\
        {\footnotesize $\langle b,c,d\rangle\leftrightarrow\langle
x,y,z\rangle$}
            \end{tabular}}

\psfrag{N}{\begin{tabular}{c}{\footnotesize $\langle
a,b,c\rangle\leftrightarrow\langle y,z,x\rangle$}\\
{\footnotesize $\langle d,b,c\rangle\leftrightarrow\langle
y,z,x\rangle$}
\end{tabular}} \psfrag{O}{\begin{tabular}{c}{\footnotesize $\langle
a,b,c\rangle\leftrightarrow\langle z,x,y\rangle$}\\
{\footnotesize $\langle c,d,b\rangle\leftrightarrow\langle
z,x,y\rangle$}
\end{tabular}} \psfrag{P}{\begin{tabular}{c}(A)\hspace{.025
in}{\footnotesize $\langle a,b,c\rangle\leftrightarrow\langle
b,c,d\rangle$}
 \end{tabular}}
\psfrag{Q}{\begin{tabular}{c}(B)\hspace{.025 in}{\footnotesize
$\langle
a,b,c\rangle\leftrightarrow\langle c,d,b\rangle$}\\
 \end{tabular}}
\psfrag{R}{\begin{tabular}{c}(C)\hspace{.025 in}{\footnotesize
$\langle a,b,c\rangle\leftrightarrow\langle d,b,c\rangle$}
 \end{tabular}}
\psfrag{L}{Layer}
        \vspace{0 in}
        \begin{center}
\epsfxsize = 3.5 in \epsfbox{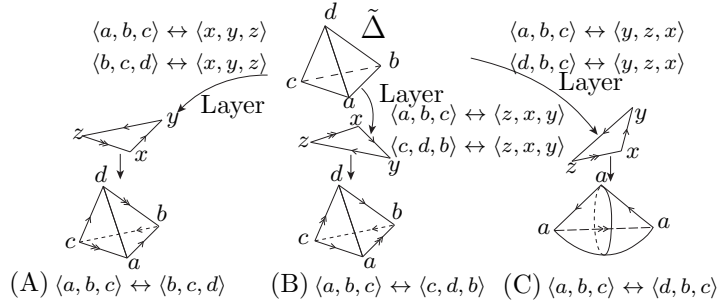}
 \caption{One-tetrahedron solid torus and creased $3$--cell (layering
  of a tetrahedron on a one-triangle M\"obius band).} \label{f-one-tet-layer-planar}
\end{center}
\end{figure}

There is another form of layering, which can be thought of as a
degenerate form of what we have just described. For example, there
are three ways to layer the back two faces of the tetrahedron
$\td{\Delta}$ onto the one-triangle M\"obius band. See Figure
\ref{f-one-tet-layer-planar}.

In parts (A) and (B) the tetrahedron is layered along the interior
(orientation reversing edge) on the one-triangle M\"obius band, the
labels and arrows give the identifications. Combinatorially, these
triangulations are the same. This triangulation of the solid torus
will be referred to as {\it the one-tetrahedron solid torus}.  In
the last case, Figure \ref{f-one-tet-layer-planar}, Part C, we show
a {\it creased $3$--cell} obtained by a single layering of a
tetrahedron along the boundary edge of the one-triangle M\"obius
band; again the labels and arrows give the identification. The
M\"obius band and the creased $3$--cell both have the homotopy type
of a solid torus; we think of each as a degenerate
layered-triangulation of the solid torus.

\begin{figure}[htbp]

            \psfrag{e}{\Large{$e$}}
            \psfrag{D}{\Large{$\td{\Delta}_t$}}
            \psfrag{f}{$\td{e}$}
            \psfrag{e}{$e$}
            \psfrag{k}{layered}\psfrag{T}{\begin{tabular}{c}
            {\small M\"obius}\\
        {\small band}\\
\end{tabular}}
             \psfrag{l}{\begin{tabular}{c}
            {\small layered}\\
        {\small solid}\\
        {\small torus}\\
            \end{tabular}}
            \psfrag{s}{\begin{tabular}{c}
            $t$\\
       layers\\
            \end{tabular}}
            \psfrag{t}{\begin{tabular}{c}
            $(t-1)$\\
            layers\\
            \end{tabular}}

        \vspace{0 in}
        \begin{center}
\epsfxsize =2.75 in \epsfbox{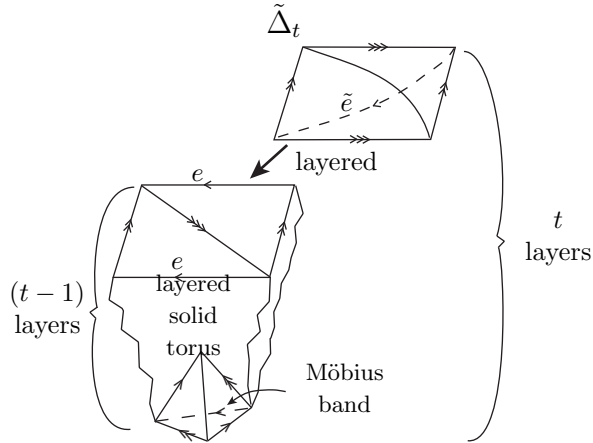}
\caption{Layered-triangulation of a solid torus.}
\label{f-layered-torus-planar-def}
\end{center}
\end{figure}

With the above notion of layering on a triangulation and starting
with the one-triangle M\"obius band,  we inductively define a
triangulation $\T_t$ of the solid torus to be a {\it
layered-triangulation of the solid torus with $t$--layers} if
\begin{enumerate}
\item[(0)] $\T_0$ is the one-triangle M\"obius band,\item[(1)] $\T_1$ is
either the one-tetrahedron solid torus or the creased $3$--cell
(both obtained by layering on the one-triangle M\"obius band), and
\item[(t)] $\T_{t} = \T_{t-1} \cup_e \td{\Delta}_t$ is a layering along
the edge $e$ of a layered-triangulation $\T_{t-1}$ of the solid
torus having $t-1$ layers, $t \geq 1$. See Figure
\ref{f-layered-torus-planar-def}.
\end{enumerate}

Note that a layered-triangulation of the solid torus with $t$ layers
has $t$ tetrahedra, $2t+1$ faces with two faces in the boundary,
$t+2$ edges with three edges in the boundary and one vertex, which
is in the boundary. It is possible by layering on the creased
$3$--cell that one does not get a solid torus but a homotopy solid
torus; we do not use such layerings.

While there is a unique one-vertex triangulation of the solid torus,
there are infinitely many ways, up to homeomorphism of the solid
torus, to place a one-vertex triangulation of the torus onto the
boundary of the solid torus. Two one-vertex triangulations,
$\T_{\bdy}$ and $\T_{\bdy}'$, on the boundary of the solid torus are
equivalent if and only if there is a homeomorphism of the solid
torus taking $\T_{\bdy}$ to $\T_{\bdy}'$. In general, if $\T_{\bdy}$
is a triangulation on the boundary of a $3$--manifold $M$,  a
triangulation $\T$ of $M$ is an extension of $\T_{\bdy}$ if $\T$
restricted to $\bdy M$ is $\T_{\bdy}$ and all the vertices of $\T$
are in $\bdy M$ (no vertices are added). We have the following
theorem from \cite{jac-rub-layered}, a version also appears in
\cite{jac-sedg-dehn}.

\begin{thm}\label{extend-layered} Suppose $\T_{\bdy}$ is a one-vertex
triangulation on the boundary of the solid torus. Then $\T_{\bdy}$
can be extended to a layered-triangulation of the solid torus; in
fact, there is a unique extension of $\T_{\bdy}$ to a minimal
layered-triangulation of the solid torus.\end{thm}

Here we use {\it minimal} with layered-triangulation to mean that of
all layered-triangulations that extend the triangulation
$\T_{\bdy}$, there is a unique one with the fewest number of
tetrahedra. We distinguish the unique minimal layered-triangulation
of the solid torus that extends the triangulation  $\T_{\bdy}$ on
its boundary by saying it is a {\it
$\T_{\bdy}$--layered-triangulation of the solid torus}. It is the
preferred extension of the equivalence class of $\T_{\bdy}$ on the
boundary of the solid torus to a layered-triangulation of the solid
torus. We do not know if the $\T_{\bdy}$--layered-triangulation is
the minimal triangulation of the solid torus extending $\T_{\bdy}$.
We conjecture that it is the minimal extension of $\T_{\bdy}$.

In \cite{jac-rub-layered} the normal surfaces in a minimal
layered-triangulations of the solid torus are characterized. In this
paper we  use that for any layered-triangulation of the solid torus,
there is a unique normal meridional disk and that in a minimal
layered-triangulation of a solid torus the only normal surface with
boundary slope meridional is the meridional disk. There are no
closed normal or almost normal surfaces in a layered-triangulation
of the solid torus.

Layered-triangulations of the solid torus are quite useful for the
construction of nice triangulations of Dehn fillings of knot- and
link-manifolds.

\begin{figure}[htbp]

            \psfrag{1}{\Large{$B_1$}}\psfrag{2}{\Large{$B_2$}}\psfrag{3}{\Large{$B_3$}}
            \psfrag{K}{\Large{$B_K$}}
            \psfrag{A}{\Large{$(M(\alpha^1,\alpha^2), \T(\alpha^1,\alpha^2))$}}
            \psfrag{M}{\Large$(M,\T)$}

            \psfrag{a}{\small$\bbb{T}(\alpha^1)$}\psfrag{b}{\small$\bbb{T}(\alpha^2)$}
            \psfrag{L}{\begin{tabular}{c}
            {\small layered}\\
        {\small solid torus}\\
\end{tabular}}

            \psfrag{c}{\small$\alpha^1$}\psfrag{d}{\small$\alpha^2$}

        \vspace{0 in}
        \begin{center}
\epsfxsize =2.5 in \epsfbox{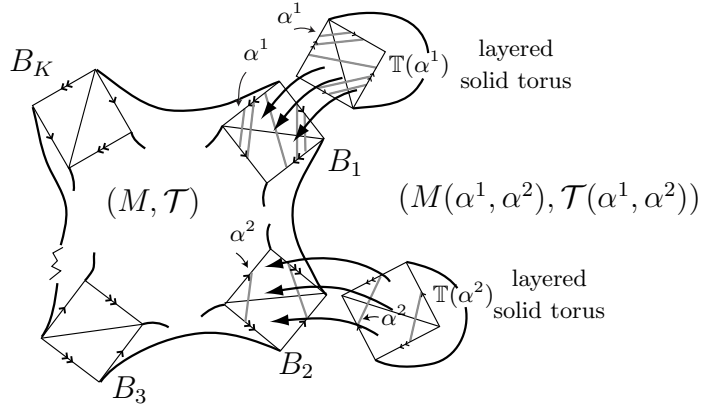} \caption{A
link-manifold with triangulation $\T$, having precisely one vertex
on each boundary component, along with layered-triangulations of the
solid torus $\bbb{T}(\alpha^1)$ and $\bbb{T}(\alpha^2)$ attached
simplicially to the boundary components $B_1$ and $B_2$. This gives
a triangulated Dehn filling $\T(\alpha^1,\alpha^2)$ of the Dehn
filling $M(\alpha^1,\alpha^2)$.} \label{f-triang-dehn-fill}
\end{center}
\end{figure}

Suppose $M$ is a link-manifold and $\T$ is a triangulation of $M$
that induces a one-vertex triangulation on each (torus) boundary
component. Suppose $\alpha$ is a slope in the component $B$ of $\bdy
M$. Let $M(\alpha)$ denote the Dehn filling of $M$ along the slope
$\alpha$ and $\bbb{T}(\alpha)$ denote the solid torus of the Dehn
filling. The triangulation $\T$ induces a triangulation $\T_B$ on $B
= M\cap\bdy\bbb{T}$. Hence, by  Theorem \ref{extend-layered}, we can
extend $\T_B$ to the $\T_B$--layered-triangulation of the solid
torus $\bbb{T}$, giving a triangulation $\T(\alpha)$ of the Dehn
filling $M(\alpha)$. We call the pair $(M(\alpha),\T(\alpha))$ a
{\it triangulated Dehn filling}; sometimes the triangulated Dehn
filling is understood by using just $M(\alpha)$ or $\T(\alpha)$. See
Figure \ref{f-triang-dehn-fill}.

If $M$ has more than one boundary component, then $M(\alpha)$ is a
link-manifold and we can consider Dehn filling it along a slope in
its boundary. Hence, if we do Dehn fillings of $M$ along slopes
$\alpha^1,\ldots,\alpha^k$ in distinct boundary components, to
arrive at the Dehn filled $3$--manifold
$M(\alpha^1,\ldots,\alpha^k)=
M(\alpha^1,\ldots,\alpha^{(k-1)})(\alpha^k)$, we can extend $\T$
using minimal layered-triangulations of the solid torus to get a
triangulated Dehn filling with triangulation
$\T(\alpha^1,\ldots,\alpha^k)=\T(\alpha^1,\ldots,\alpha^{(k-1)})(\alpha^k)$.

Notice that for triangulated Dehn fillings of link-manifolds, the
triangulation is $\T$ on $M$ and is a minimal
$\T_{B_i}$--layered-triangulation on each solid torus
$\bbb{T}(\alpha^i)$, where $B_i$ is the torus boundary component
containing the slope $\alpha^i$ and $\T_{B_i}$ is the triangulation
on $B_i$ induced by $\T$.

Now, if $\widehat{G}$ is a normal surface in
$M(\alpha^1,\ldots,\alpha^k)$ with a triangulated Dehn filling
$\T(\alpha^1,\ldots,\alpha^k)$, then $\widehat{G}$ meets $M$ and
each of the solid tori $\bbb{T}(\alpha^i)$ in normal surfaces. This
is an aspect of these triangulations, along with the
characterization of normal surfaces in a minimal
layered-triangulation of a solid torus, that is helpful in
understanding normal and almost normal surfaces in Dehn fillings. If
the components of $\widehat{G}$ in each of the solid tori
$\bbb{T}(\alpha^i)$ are disks and $G$ is the normal surface in which
$\widehat{G}$ meets $M$, we say $G$ ``caps off" and write
$G(\alpha^1,\ldots,\alpha^k)$ for $\widehat{G}$. Of course, there
are normal surfaces in the triangulated Dehn filling
$(M(\alpha^1,\ldots,\alpha^k),\T(\alpha^1,\ldots,\alpha^k))$ that
meet the solid tori in normal surfaces other than the meridional
disks; but the ``capped-off" normal surfaces are special. We have
the following.

\begin{lem} \label{re-write}Suppose
$(M(\alpha^1,\ldots,\alpha^k),\T(\alpha^1,\ldots,\alpha^k))$ is a
triangulated Dehn filling. If the ``capped-off" normal surface
$P(\alpha^1,\ldots,\alpha^k)=\sum n_q\widehat{G}_q$ is a geometric
sum of normal surfaces in
$(M(\alpha^1,\ldots,\alpha^k),\T(\alpha^1,\ldots,\alpha^k))$, then
$\widehat{G}_q=G_q(\alpha^1,\ldots,\alpha^k)$ is a ``capped off"
surface for each $q, 1\le q\le k$, and $P=\sum n_qG_q$.\end{lem}
\begin{proof} If $\bbb{T}(\alpha^i)$ is one of the layered solid tori
in the triangulated Dehn filling, then $P(\alpha^1,\ldots,\alpha^k)$
meets $\bbb{T}(\alpha^i)$ in a family of meridional disks;
furthermore, the components of intersection of each $\widehat{G}_q$
with $\bbb{T}(\alpha^i)$ are normal surfaces and have geometric sum
the collection of meridional disks in which
$P(\alpha^1,\ldots,\alpha^k)$ meets $\bbb{T}(\alpha^i)$. It follows
from Theorem \ref{slopes} that the components of each
$\widehat{G}_q$ meet $\bdy\bbb{T}(\alpha^i)$ in the meridional slope
$\alpha^i$. However, by the classification of normal surfaces in a
minimal layered-triangulation of a solid torus, the only such normal
surface in $\bbb{T}(\alpha^i)$ is the unique normal meridian disk.
Thus each $\widehat{G}_q$ ``caps off" and if $\widehat{G}_q =
G_q(\alpha^1,\ldots,\alpha^k)$, then we have $P = \sum
n_qG_q$.\end{proof}

In the situation above, we say we can {\it re-write} $P$ as the sum
$\sum n_qG_q$. This is particularly significant when we have $P$
written as a sum $P=\sum m_jF_j$, where the $F_j$ are fundamental in
$(M,\T)$. Then we re-write $P =\sum n_qG_q$ where
$G_q(\alpha^1,\ldots,\alpha^k)$ is fundamental in
$(M(\alpha^1,\ldots,\alpha^k),\T(\alpha^1,\ldots,\alpha^k))$. In
re-writing $P$ in this way, we have each of the summands $G_q$
meeting the boundary component with slope $\alpha^i$ in the slope
$\alpha^i$, whereas, there is no way of knowing how the various
$F_j$ meet the boundary component containing the slope $\alpha^i$ in
a link-manifold with multiple boundary components.

\subsection{Basic algorithms.} In this subsection, we organize some
of the algorithms we will be using. Note that we often state the
existence of an algorithm about surfaces in a $3$--manifold by
stating that there is an algorithm that will decide; and if the
answer is yes, then the algorithm will construct the desired
surface. In many case, however, this is done by showing that the
answer is yes if and only if there is a {\it fundamental surface}
that is of the type we seek. And, of course, fundamental surfaces
can be constructed. An interesting aspect of the main results of
this paper is that the algorithms we develop later do not
necessarily find an answer among the fundamental surfaces; to get an
answer we have to look, in some cases, rather far afield. We shall,
generally, state the conclusions of our results in the stronger
terms of finding solutions among the fundamental surfaces, if,
indeed, that is where a solution can be found. Many of the results
we use here are improved in \cite{jac-tol} and \cite{jac-ree},
showing that desired solutions are already at the vertices of
projective solution space if they exist at all. Also, in some cases,
a conclusion can be made for any triangulation; in other cases, we
must first modify the given triangulation to a triangulation more
suitable to a solution of the problem.

The first result we give is generally attributed to W. Haken
\cite{haken1}. We note that the proof attributed to Haken is for a
manifold given via a handle decomposition where it is known that the
manifold is orientable and irreducible. Such a proof, learned from
the methods of \cite{sch}, is given in \cite{jac-oer}. In G.
Hemion's book, \cite{hemion-book}, an argument is given in the case
of triangulations; however, it parallels the argument for handle
decompositions, leaving it with a gap. The gap is eliminated in
\cite{jac-ree} where the following result is proved. We remark that
the argument in this generality does involve many details; we can
obtain the result much more easily by first constructing a prime
decomposition of the manifold from which we will find an essential
disk should one exist (see Theorem \ref{prime-decomp} below).

\begin{lem}\label{find-disk} Suppose $M$ is a compact $3$--manifold
and $B$ is a component of $\bdy M$. $M$ contains a properly
embedded, essential disk with boundary in $B$ if and only if for any
triangulation $\T$ of $M$, there is an essential, normal disk with
boundary in $B$ among the fundamental solutions for
$(M,\T)$.\end{lem}

 If we are seeking to know the existence of
essential surfaces, then if they exist, they will either be part of
a prime decomposition (spheres, projective planes) or exist only if
they are in prime factors. In fact, later we reduce our main problem
to a problem where it is known that the given manifold is
irreducible and $\bdy$--irreducible. However, in the steps of an
algorithm, it often is necessary to alter the given knot or
link-manifold to one obtained by a Dehn filling of the given
manifold. It is well known that upon Dehn filling one might lose
some of the nice features of the original given manifold, such as
irreducibility and $\bdy$--irreducibility.

A version of the following lemma appears within the proof of Lemma
5.11 of \cite{jac-sedg-dehn}; however, that version assumes it is
known that the given manifold is irreducible and
$\bdy$--irreducible.

\begin{lem}\label{find-annulus} Suppose $M$ is a compact $3$--manifold with nonempty
boundary.  $M$ contains a properly embedded, essential annulus
having its boundary  in distinct boundary components $B$ and $B'$ of
$M$ if and only if for any triangulation $\T$ of $M$, there is an
embedded, essential, fundamental normal annulus having its boundary
components in $B$ and $B'$.
\end{lem}

\begin{rem} Suppose $M$ is link-manifold and $B$ is a component of $\bdy
M$. If there are properly embedded, essential annuli $A$ and $A'$ in
$M$ each having a boundary component in $B$, then either the
component(s) of the boundary of $A$ and $A'$ in $B$ have the same
slope or $B$ is in the boundary of a prime factor of $M$ that is an
$I$-bundle; i.e., $S^1\times S^1\times I$ or the twisted $I$--bundle
over the Klein bottle.\end{rem}

A general version of an algorithm to decide if a given normal
surface is an essential surface appears in \cite{jac-tol}; the
following version follows from early work of Haken \cite{haken1}
with some of the material from \cite{jac-tol} to determine if a
normal surface in a $3$--manifold $M$ is parallel into $\bdy M$.
Recall in a link-manifold $M$ a properly embedded surface is
essential if and only if it is incompressible and is not an annulus
or torus parallel into $\bdy M$. One does not need separately to
check $\bdy$--irreducibility.

\begin{thm}\label{decide-essential} Given a link-manifold
$M$ and an embedded, normal surface $F$ in $M$, there is an
algorithm to decide if $F$ is an essential surface in $M$.\end{thm}

We also have another very useful result that we can use to conclude
that a surface is essential. It was first established in
\cite{jac-oer} but using handle-decompositions; later it was redone
for triangulations in \cite{bart-sch} and in the form we use in
\cite{jac-tol}.

\begin{thm}\label{carrier-essential} Suppose $M$ is an irreducible, $\bdy$--irreducible $3$--manifold
and $\T$ is a triangulation of $M$. If $F$ is an embedded, essential
normal surface in $M$ and is least weight in its isotopy class, then
every normal surface with projective class in the carrier of $F$,
  $\C(F)$, is embedded and essential in $M$.\end{thm}

Later, Tollefson \cite{tol-iso} showed that all normal surfaces with
projective class in $\C(F)$ are also least weight in their isotopy
class.

Finally, given a normal surface, there are various ways to determine
its Euler characteristic, its connectivity, its orientability class
and the number of its boundary components. Hence, its genus (if
orientable, the number of handles; and if non-orientable, the number
of cross caps) also can be determined.

\section{Average Length Estimates} In this section we give a tool
which is very useful in working with decision problems and
algorithms, especially those related to Dehn fillings. We call it an
average length estimate; it is used extensively in
\cite{jac-sedg-dehn}. If $M$ is a $3$--manifold with triangulation
$\T$ and $F$ is a properly embedded surface in $M$ and is in general
position with respect to $\T^{(2)}$, then we have $L(\bdy F)$
defined. If $b$ is the number of components of $\bdy F$, then we set
$\lambda_{av} = L(\bdy F)/b$ and say $\lambda_{av}$ is the {\it
average length (of the components) of $\bdy F$}. Under various
conditions placed on the topology of a manifold $M$, the average
length estimate says that for any triangulation $\T$ of $M$, there
is a constant $C$, depending only on $M$ and $\T$, so that all
properly embedded, essential surfaces of bounded genus have the
average length of their boundary bounded by $C$. Hence, depending
only on the manifold and a given triangulation, there is a number so
that all essential surfaces of bounded genus must have a short
boundary component. We give several variants here and add another
useful variant within the proof of Theorem \ref{find-planar-bdd}.

If a $3$--manifold has an essential annulus, then often it is
possible by Dehn twisting about such an annulus to obtain surfaces
with all boundary components being arbitrarily long; the surfaces
are homeomorphic and so all have a fixed genus. Hence, there are no
preassigned values for short boundaries. There also are examples of
families of surfaces of fixed genus in link-manifolds, where for any
value there is a surface in the family having some boundary
component of length larger than this value; this can happen in a
link-manifold with or without having an essential annulus between
distinct boundary components. However, if there are no such
essential annuli, then each such surface must also have short
boundary components.  We say the $3$--manifold $M$ is {\it
anannular} if there are no properly embedded, essential annuli in
$M$.

\begin{prop}\label{ALE} Suppose $M$ is a compact, irreducible,
$\bdy$--irreducible and anannular $3$--manifold with triangulation
$\T$. There is a constant $C = C(M,\T)$, depending only on $M$ and
$\T$, so that if $F$ is an embedded, essential, normal surface that
is least weight in its isotopy class and $\lambda_{av}$ is the
average length of the components of $\bdy F$, then $\lambda_{av}\le
C(2g+1),$ where $g$ is the genus of $F$.\end{prop}
\begin{proof} Since $F$ is essential and least weight in its isotopy
class, we have from \cite{jac-oer} and \cite{jac-tol} that any
normal surface that projects into $\C(F)$ is essential. In
particular, from the hypotheses on $M$, no normal surface in $\C(F)$
is a $2$--sphere, projective plane, annulus or M\"obius band. Hence,
the normal surface $F$ can be written as a sum $F = \sum n_iF_i
+\sum m_jK_j $, where each $F_i$ is fundamental and essential and
each $K_j$ is an essential torus or Klein bottle.

Set
$$C = \frac{L(\bdy F_i)}{-\chi(F_i)},$$ where $F_i$ is a fundamental
surface for $(M,\T)$ and $\chi(F_i)< 0$.

Then we have  $$b\lambda_{av} = L(\bdy F) = \sum n_i L(\bdy F_i)\le
C\sum n_i(-\chi(F_i)) = C(-\chi(F))=C(2g-2+b),$$ where $b$ is the
number of components of $\bdy F$. It follows that $\lambda_{av} \le
C(2g+1).$ \end{proof}

We have noted that the possibility of Dehn twisting a surface $F$
about a properly embedded annulus in the manifold can lengthen the
boundary of $F$ but does not change its genus. However, after
passing to normal surfaces,  this feature can only happen if such an
annulus is fundamental and also has the same quadrilateral type as
the surface $F$. Also, Dehn twisting a surface $F$ about an annulus
having its boundary in a single boundary component that is a torus
does not change the slope of $\bdy F$. We have two variants to
Proposition \ref{ALE} given as the next proposition and its
corollary. In addition, we point out that while Euler characteristic
is additive under geometric sum, genus, in general, is not. The
problem, of course, is that if a normal surface is a geometric sum
of other normal surfaces, then it is not always true that the number
of boundaries of the sum is the sum of the number of boundaries of
the summands. But in many cases Theorem \ref{slopes} does allow this
to happen and thus gives us results when the manifold has tori
boundary that we would not get otherwise. This is quite evident in
the next two results.

\begin{prop}\label{ALE-link1} Suppose $M$ is a
link-manifold with no embedded annuli having essential boundary
curves in distinct components of $\bdy M$. Furthermore, suppose $\T$
is a $0$--efficient triangulation of $M$. Then there is a constant
$C=C(M,\T)$, depending only on $M$ and $\T$, so that if $F$ is an
embedded normal surface in $M$ with no trivial boundary curves and
$\lambda_{av}$ is the average length of the components of $\bdy F$,
then $\lambda_{av}\le 2C(g+1),$ where $g$ is the genus of
$F$.\end{prop}\begin{proof} If $F$ is an embedded normal surface
with no trivial boundary curves, then we can write $$F = \sum
l_iF_i+\sum m_jK_j+\sum n_kA_k,$$ where each summand is fundamental,
$\chi(F_i)<0, K_j$ is either a torus or Klein bottle and $A_k$ is
either a M\"obius band or an annulus with both its boundary curves
in the same component of $\bdy M$. We let $\abs{A_k}$ denote the
number of boundary components of $A_k$ and set $$C =
max\left\{\frac{L(\bdy F_1)}{-\chi(F_1)},\ldots,\frac{L(\bdy
F_I)}{-\chi(F_I)}, \frac {L(\bdy A_1)}{\abs{A_1}},\ldots,\frac
{L(\bdy A_K)}{\abs{A_K}}\right\}.$$ Let $F' = \sum l_iF_i$; then
$\chi(F) = \chi(F')$.

Now, by Theorem \ref{slopes}, we have that if an annulus or M\"obius
band summand and $F'$ meet the same boundary torus of $M$, then
their boundaries have the same slope in this boundary torus. Hence,
if $b_{F'}$ is the number of boundary components of $\bdy F'$, $b_A
= \sum n_k\abs{A_k}$, is the number of boundary components of $\sum
n_kA_k$ and $b$ is the number of boundary components of $F$, then $b
= b_{F'}+b_A$. It follows that
$$b\lambda_{av} = L(\bdy F) = \sum l_iL(\bdy F_i) + \sum
n_kL(\bdy A_k)\le$$ $$\le C\left(\sum l_i(-\chi(F_i))+\sum
n_k\abs{\bdy A_k}\right)=C(-\chi(F') + b_A).$$ However, $\chi(F) =
\chi(F')$; hence,  $b\lambda_{av}\le C(-\chi(F)+b_A)= C(2g-2
+b)+Cb_A$. Therefore, $$\lambda_{av}\le  C(\frac{(2g-2)}{b}
+1)+C(\frac{b_A}{b})\le C(2g+2).$$\end{proof}

Suppose $M$ is a $3$--manifold and $\T$ is a triangulation of $M$.
We shall say a finite collection of embedded normal surfaces
$\{H_1,\ldots,H_m\}$ is a {\it spanning collection} for the embedded
normal surfaces in $(M,\T)$ if and only if for any embedded normal
surface $F$ in $M$, we have $F = \sum n_iH_i$, where $n_i$ is a
nonnegative integer. The fundamental surfaces are a spanning
collection.

The following lemma is a corollary of the proof of the Proposition
\ref{ALE-link1}. In this corollary, we  substitute a spanning
collection for the fundamental surfaces; hence, the constant $C$
becomes dependent on the spanning collection (a fundamental
collection is unique; whereas, there are possibly many choices of
distinct spanning collections).

\begin{cor}\label{ALE-link2} Suppose $M$ is
link-manifold, $\T$ is a minimal-vertex triangulation of $M$ and
$\mathcal{H}$ is a spanning collection of normal surfaces. There is
a constant $C=C(M,\T, \mathcal{H})$, depending only on $M$, $\T$ and
$\mathcal{H}$, so that if $F$ is an embedded normal surface in $M$
with no trivial boundary curves and $F$ can be written as a
geometric sum $F = \sum l_iH_i$, where $H_i \in\mathcal{H}$  and
either $\chi(H_i)<0$ or $H_i$ is a torus or Klein bottle, then
$$\lambda_{av}\le C(2g+1),$$ where $\lambda_{av}$ is the average length
of the components of $\bdy F$ and $g$ is the genus of $F$.\end{cor}

\section{finding planar surfaces}\label{planar}

In this section we give an algorithm to decide if a given
link-manifold contains a properly embedded, essential, planar
surface. An important aspect of this algorithm is, unlike those in
Section 2, we do not show that if there is such a planar surface,
then there is one among the fundamental surfaces of the given
triangulation. Indeed, we must construct a family that, while still
finite, goes beyond the fundamental surfaces of the given
triangulation. The background material and other algorithms we use
in the proof of this result typically assume the given manifold is
known to be irreducible and $\bdy$--irreducible. So, we begin with
an algorithm that transforms the problem for the given manifold to a
possibly distinct but constructible manifold that is known to be
irreducible and $\bdy$--irreducible.

Beginning in the most general situation, we are given a
$3$--manifold via a triangulation. We can easily check that it is a
link-manifold; however, it very well may  be $\bdy$--reducible or
reducible. Recall that if a link-manifold $M$ is $\bdy$--reducible
and is not a solid torus, then it also is reducible and can be
written as a connected sum, $M = (D^2\times S^1)\#M'$, of a solid
torus and a $3$--manifold $M'$. There is an algorithm, due to W.
Haken \cite{haken1}, which in Lemma \ref{find-disk} we adapted to
the generality we are using here, that will decide if a given
$3$--manifold is $\bdy$--reducible. If it is, the algorithm will
construct an essential disk. Hence, we could begin by running this
algorithm. Having run this algorithm, if the given link-manifold is
$\bdy$--reducible, we have found a properly embedded, essential
planar surface and we are done. However, if the manifold is not
$\bdy$--reducible, it may still be reducible. Thus it would still be
necessary to undertake the construction of a prime decomposition.
So, we do this first, including, with little extra work, an
algorithm that will modify a given triangulation of a manifold that
is known to be irreducible and $\bdy$--irreducible to a
triangulation of the manifold that is $0$-efficient. We note,
however, there is a version of the Haken algorithm within these
algorithms. Also, we do not really need that triangulations are
$0$--efficient; but rather use very strongly that the given
manifolds are irreducible and $\bdy$--irreducible and are given by
minimal-vertex triangulations.

If $Q$ is a compact $3$--manifold, we use the notation $p(Q)$ to
denote $Q\#\cdots\#Q$, where there are $p\ge 0$ copies of $Q$.

\begin{thm}\cite{jac-rub0}\label{prime-decomp} Given a link-manifold $M$ via a triangulation
$\T$, there is an algorithm to construct a prime decomposition
$$M=p(S^2\times S^1)\# q(\rppp)\# r(D^2\times S^1)\#
M_1\#\cdots\#M_n,$$ where $p,q$ and $r$ are nonnegative integers and
each $M_i$ is given by a $0$--efficient triangulation $\T_i,  i
=1,\ldots,n$, respectively; furthermore, $\sum_{i=1}^n Card(\T_i)\le
Card(\T)$.\end{thm}
\begin{proof}  First, from
\cite{jac-ree} and Proposition 5.7 of \cite{jac-rub0}, there is an
algorithm to decide if there are non vertex-linking normal
$2$--spheres. If there is one, the algorithm constructs one. We then
proceed as in the proof of Theorem 5.9 of \cite{jac-rub0} to crush
the triangulation along a suitable non vertex-linking normal
$2$--sphere. In the case of a link-manifold, if there is a non
vertex-linking normal $2$--sphere, there is one along which we can
crush. The process repeats as long as there is a non vertex-linking
normal $2$--sphere. After each crushing, we reduce the number of
tetrahedra we had before crushing, thus the process must stop in a
finite number of steps at which point the algorithm has managed to
find any factors that are $S^2\times S^1$ or $\rppp$ and we have a
connected sum decomposition
$$M=p(S^2\times S^1)\# q(\rppp)\#
M'_1\#\cdots\#M'_{n'},$$ where each factor $M_i'$ is given by a
triangulation in which the only normal $2$--spheres are
vertex-linking. We note that the triangulations of the closed
factors in this decomposition are $0$--efficient; and by using the
$3$--sphere recognition algorithm \cite{rubin1, tho} (also, see
\cite{jac-rub0}, Theorem 5.11), we may assume no factor is $S^3$.

If the only normal $2$--spheres are vertex-linking, then we consider
the factors $M_i'$ in the above connected sum decomposition that
have nonempty boundary. For any such factor, we determine if there
are any non vertex-linking normal disks. Again, from \cite{jac-ree}
and \cite{jac-rub0} there is an algorithm to decide if there are any
non vertex-linking normal disks. If there is one, then the algorithm
constructs one.  However, here there is a slight twist to the
argument in \cite{jac-rub0}; namely, we may have a non separating
essential disk. Since each factor is irreducible, we conclude that a
factor with a non separating disk is a solid torus and contributes a
factor in the prime decomposition of the form $D^2\times S^1$. If
the disk is separating, we proceed as in the proof of Theorem 5.17
of \cite{jac-rub0} to crush the triangulation along a suitable
separating, non vertex-linking normal disk. The process repeats as
long as there is a non vertex-linking normal disk. Again, since at
each crushing we must reduce the number of tetrahedra we had before
crushing, the process must stop in a finite number of steps.

Upon having no non vertex-linking normal disk, we have the desired
decomposition, where the factors $M_1,\ldots,M_n$ are given by
$0$--efficient triangulations $\T_1,\ldots,\T_n$, respectively.
Furthermore, the total number of tetrahedra in the triangulations
$\T_1,\ldots,\T_n$ is no larger than the number we started with in
$\T$ and is equal if and only if $\T$ is itself
$0$--efficient.\end{proof}

Knowing that we can construct a prime decomposition of a given
manifold, we make the following observation, a proof of which is
easily derived from classical $3$--manifold ``cut-and-paste"
methods.

\begin{prop}\label{planar-prime-factor} The $3$--manifold $M$ contains a properly embedded,
essential, planar surface if and only if one of the prime factors of
$M$ contains a properly embedded, essential, planar surface.
Moreover, if $B$ is a component of $\bdy M$, then $M$ contains an
embedded punctured disk with boundary in $B$ and punctures in $\bdy
M\setminus B$ if and only if the prime factor, say $M_B$ of $M$
containing $B$ has an embedded punctured-disk with boundary in $B$
and punctures in $\bdy M_B\setminus B$.
\end{prop}

If $M$ is a link-manifold, then  the decomposition given in Theorem
\ref{prime-decomp} necessarily has prime factors that are known to
be solid tori or are given by $0$--efficient triangulations and thus
are known to be irreducible and $\bdy$--irreducible link-manifolds.
We ignore any factors that are closed.

We now have a special case of the main theorem of this section;
which also is needed in the proof of the main theorem. Notice that
in Lemma \ref{all-in-one-bdry}, we do not require that the
link-manifold $M$ be irreducible or $\bdy$--irreducible and we do
not require that we consider only essential, planar surfaces. We do,
however, require that all of the boundary components of our planar
surfaces are essential curves in a single component $B$ of $\bdy M$.
There is, of course, something subtle here, as it is quite easy to
find an embedded, planar, normal surface with all of its boundary
components essential curves in a component $B$ of $\bdy M$. For
example, consider the surface one obtains by taking the frontier of
a small regular neighborhood of an edge in the component $B$ of
$\bdy M$. Often, this surface is normal and fundamental. If it is
not normal, it can be shrunk, using a barrier surface argument
\cite{jac-rub0}, it either normalizes to an embedded, normal, planar
surface (annulus) with its boundary essential curves in $B$ or it
follows that $M$ is a  solid torus. The point, which will give us a
nontrivial conclusion in Lemma \ref{all-in-one-bdry}, is that the
existence of a certain special planar, normal surface, say $P$,
guarantees not only the existence of a similar planar, fundamental
normal surface but one that also projects into $\C(P)$, the carrier
of $P$. Following the statement and proof in this special case, we
give two corollaries. The second corollary is our main theorem of
this section for knot-manifolds.

\begin{lem}\label{all-in-one-bdry} Suppose $M$ is a link-manifold  with a
minimal-vertex triangulation $\T$ and $B$ is a component of $\bdy
M$. If there is an embedded, planar, normal surface $P$ with all its
boundary essential curves in $B$ and if $P$ can be written as a sum
of fundamental surfaces, none of which is a $2$--sphere or
projective plane, then there is an embedded, planar, normal surface
with all its boundary essential curves in $B$ that projects into a
fundamental class in $\C(P)$, the carrier of $P$.\end{lem}
\begin{proof} Suppose there is an embedded, planar, normal surface
$P$ with all of its boundary essential curves in $B$ and
 $$P = \sum n_iF_i+\sum n_j'F_j',$$ where $F_i$ is
fundamental, orientable and not a $2$--sphere and  $F_j'$ is
fundamental, non-orientable and not a projective plane. Considering
Euler characteristics, we have for $b$ the number of boundary
components of $P$, $b_i, b_j'$ the number of boundary components of
$F_i, F_j'$, respectively, $$b-2 =-\chi(P)=\sum n_i(-\chi(F_i))+\sum
n'_j(-\chi (F_j')) =$$ $$=2\left(\sum n_i(g_i-1)+\sum
n_j'(\frac{c_j}{2}-1)\right)+\sum n_ib_i+\sum n_j'b_j',$$ where
$g_i$ is the genus of $F_i$ and $c_j$ is the number of cross caps in
$F_j'$. Since all the boundary curves of $P$ are in $B$, we have all
the boundary curves of the surfaces $F_i$ and $F_j'$ also in $B$. We
have all the surfaces $F_i, F_j'$ and $P$ satisfying the same
quadrilateral conditions; hence, by Theorem \ref{slopes} the slopes
of the boundary curves are complementary. But since the components
of $\bdy P$ are essential curves in $B$, all the surfaces actually
have the same boundary slope and this slope is the slope of the
boundary of $P$. It follows that $b = \sum n_ib_i+\sum n_j'b_j'$ and
therefore
$$-1 = \sum n_i(g_i-1)+\sum n_j'(\frac{c_j}{2}-1).$$

We conclude some $g_i =0$ or some $c_j =1$. If $g_i = 0$, then $F_i$
is a fundamental planar surface (no $F_i$ is a $2$--sphere); if $c_j
= 1$, then $F_j'$ is a fundamental  M\"obius band (no $F_j'$ is a
projective plane). All summands in the geometric sum of $P$ must
project into $\C(P)$.\end{proof}

\begin{cor} Given a link-manifold  $M$, there is an algorithm to decide if $M$
contains a properly embedded, essential, planar surface with all its
boundary in a given component $B$ of $\bdy M$. If there is one, the
algorithm will construct one.\end{cor}

\begin{proof} We are given the link-manifold $M$ via a triangulation
$\T$. By Theorem \ref{prime-decomp}, we construct a prime
decomposition $$P = p(S^2\times S^1)\#q(\rppp)\#r(D^2\times
S^1)\#M_1\#\cdots\#M_n,$$ where each $M_i$ has a $0$--efficient
triangulation. If $r\ne 0$ and $B$ is a component of a solid torus
factor, $D^2\times I$, then we have an essential disk with boundary
in $B$ and the algorithm that constructs the prime decomposition
will construct such an essential disk. So, we shall assume $B$ is a
component of some $M_i$, say $M_1$.

By Lemma \ref{planar-prime-factor}, $M$ contains a properly
embedded, essential, planar surface with all its boundary in  $B$ if
and only if  $M_1$ contains a properly embedded, essential, planar
surface with all its boundary in  $B$. Furthermore, if $M_1$
contains such a surface then it contains a normal one; and from
Theorem \ref{carrier-essential}, if $P$ is the least weight properly
embedded, essential, planar normal surface with all its boundary in
$B$, then every surface in $\C(P)$ is essential. Since $M_1$ has a
$0$--efficient triangulation, $P$ can be written as a sum of
fundamental normal surfaces, each must project into $\C(P)$, and
none can be a $2$--sphere or a projective plane. Thus by Lemma
\ref{all-in-one-bdry}, there is a planar surface projectively
equivalent to the projection of a fundamental surface into $\C(P)$.
Every surface that projects into $\C(P)$ is essential.
\end{proof}

\begin{cor}\label{knot-planar} Given a knot-manifold there is an algorithm to decide if
it contains a properly embedded, essential, planar surface; if it
does, the algorithm will construct one.\end{cor}

\begin{rem} If in either of the situations given in Corollaries 4.4 or 4.5
we know the given manifold is irreducible and if it is given by a
minimal-vertex triangulation $\T$, then it contains a properly
embedded, essential, planar surface if and only if there is an
embedded, essential planar, normal surface that is fundamental in
$(X,\T)$.
\end{rem}

\begin{thm}\label{find-planar} Given a link-manifold there is an algorithm to decide if
it contains a properly embedded, essential, planar surface; if it
does, the algorithm will construct one.\end{thm}

\begin{proof} From the above, we may assume we are given a
link-manifold $M$ via a $0$--efficient triangulation  $\T$.

If there is a properly embedded, essential, planar surface in $M$,
then for any triangulation of $M$ there is a normal such surface
and, therefore, an  embedded, essential, planar, normal surface that
is least weight in its isotopy class.

The plan of the proof is to study what the situation is under the
assumption that there is such a normal surface, say $P$, in $M$. We
have from Corollary \ref{knot-planar}, that for this situation in a
knot-manifold, there is an essential, planar surface among the
fundamental surfaces for $(X,\T)$. This does not seem to be
necessarily true for link-manifolds. So, we start by constructing
the fundamental surfaces of $M$. If we find an essential, planar
normal surface among these surfaces, then we are done; if we do not
find an essential, planar surface among the fundamental surfaces of
$M$, we work our way via Dehn fillings toward a knot-manifold where
we then, hopefully, can solve the problem. The first issue is to
determine those Dehn fillings we should use. For this we use the
method of average length estimates, which depend only on the
triangulation when considering fixed genus (planar) surfaces.
However, we immediately run into the classical issues with Dehn
fillings; namely, after Dehn filling we may loose, irreducibility,
$\bdy$--irreducibility and that the surfaces we are interested in
are essential. To handle these problems, we use triangulated Dehn
fillings. This avoids ever having to re-triangulate $M$ and enables
us to understand the normal surfaces in the Dehn filled manifold
relative to normal surfaces in the manifold before Dehn filling.  In
this way, if the link-manifold has $n$ tori in its boundary, then in
a succession of no more than $n$ steps, we construct at most $n$
finite collections of normal surfaces in $(X,\T)$, showing that
there is a properly embedded, essential, planar surface in $M$ if
and only if there is an essential, planar, normal surface in the
collection of surfaces we construct.

The notation is a bit tricky. At each new step, we construct a
finite family of triangulated Dehn fillings for each member of  a
previously constructed finite family of triangulated Dehn fillings.
Then for each new construction, we compute the fundamental normal
surfaces in the new triangulated Dehn filling and select from these
a subcollection, the members of which have a particularly nice
decomposition in terms of that  Dehn filling. These decompositions
provide normal surfaces in $M$ that become part of our desired
collection of normal surfaces in $(M,\T)$. The triangulation $\T$ of
$M$ in all these triangulated Dehn fillings remains constant.

We have provide Figure \ref{f-rewrite-planar} as an example having
at most three steps.

\begin{figure}[htbp]

            \psfrag{1}{\Large{$B_1$}}\psfrag{2}{\Large{$B_2$}}\psfrag{3}{\Large{$B_3$}}
            \psfrag{K}{\Large{$B_K$}}\psfrag{4}{\footnotesize$H^0_1$}\psfrag{5}{\footnotesize$H^0_3$}
            \psfrag{6}{\footnotesize$H^1_3$}
            \psfrag{7}{\footnotesize$H^1_7$}\psfrag{8}{\footnotesize$H^2_1$}\psfrag{9}{\footnotesize$H^2_3$}
            \psfrag{A}{$(M(\alpha^1,\alpha^2), \T(\alpha^1,\alpha^2))$}
            \psfrag{B}{$(M(\alpha^1), \T(\alpha^1))$}
            \psfrag{M}{$(M,\T)$}\psfrag{P}{$\boldsymbol{P}$}
\psfrag{p}{\footnotesize{$P = \sum \ell_i
H^0_i$}}\psfrag{q}{\small$P=\sum\ell^1_j H^1_j$}
\psfrag{r}{\small$P=\sum\ell^2_k H^2_k$}
            \psfrag{a}{\small$\bbb{T}(\alpha^1)$}\psfrag{b}{\small$\bbb{T}(\alpha^2)$}
             \psfrag{k}{\small$\alpha^K$}
             \psfrag{i}{\begin{tabular}{c}
            {\footnotesize $P(\alpha^1) = \sum
\ell^1_j\widehat{H}^1_j$}\\
        {\footnotesize $P(\alpha^1) = \sum
\ell^1_j{H}^1_j(\alpha^1)$}\\
\end{tabular}}
\psfrag{S}{\LARGE$\S(M,\T)$} \psfrag{j}{\begin{tabular}{c}
            {\footnotesize $P(\alpha^1,\alpha^2) = \sum
\ell^2_k\widehat{H}^2_k$}\\
        {\footnotesize $P(\alpha^1,\alpha^2) = \sum
\ell^2_k{H}^2_k(\alpha^1,\alpha^2$)}\\
\end{tabular}}
            \psfrag{L}{\begin{tabular}{c}
            {\footnotesize layered}\\
        {\footnotesize solid torus}\\
\end{tabular}}

            \psfrag{c}{\small$\alpha^1$}\psfrag{d}{\small$\alpha^2$}

        \vspace{0 in}
        \begin{center}
\epsfxsize =5 in \epsfbox{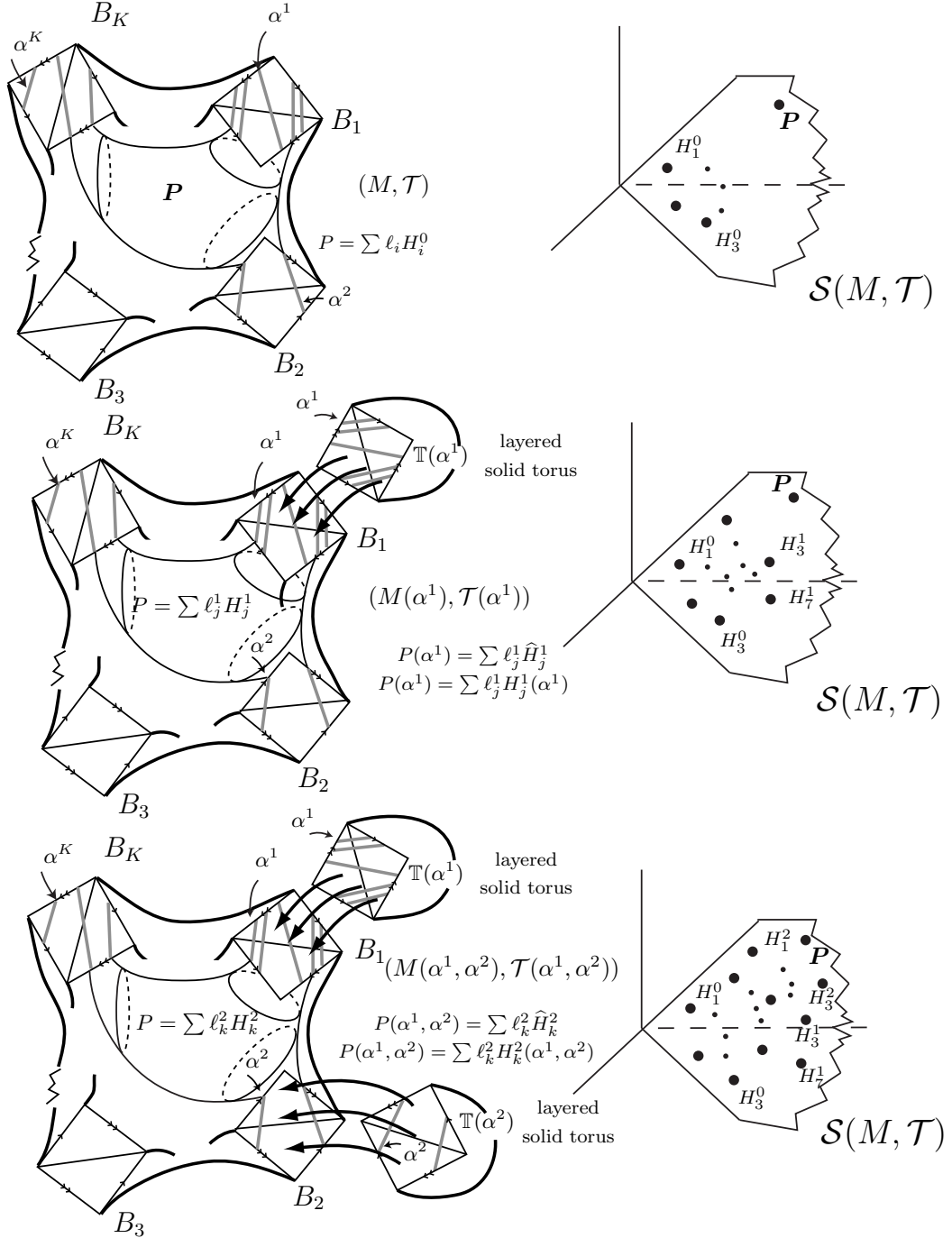} \caption{First the
planar surface $P$ is written as a sum of fundamentals $H^0_i$ in
$(M,\T)$; then we re-write $P$ as a sum $P=\sum\ell^1_j H^1_j$,
where the $H^1_j$ are determined by fundamentals in
$(M(\alpha^1),\T(\alpha^1))$; and finally we re-write $P$ as a sum
$P=\sum\ell^2_k H^2_k$, where the $H^2_k$ are determined by
fundamentals in $(M(\alpha^1,\alpha^2),\T(\alpha^1,\alpha^2))$.
There is an essential planar surface in $M$ iff there is one in one
of the finite  collections $\{H^0_i\}, \{H^1_j\}$ or $\{H^2_k\}$.}
\label{f-rewrite-planar}
\end{center}
\end{figure}

To this end, we suppose $P$ is an embedded, essential, planar,
normal surface and is least weight in its isotopy class. Let $\C(P)$
denote the carrier of $P$. Then by Theorem \ref{carrier-essential}
(see \cite{jac-oer} and \cite{jac-tol}) any normal surface in
$(M,\T)$ that is projectively equivalent to a surface in $\C(P)$ is
essential in $M$.

Consider the fundamental surfaces in $(M,\T)$ with projective class
in $\C(P)$. Such a surface is essential in $M$ and can not be a
disk, $\rpp$ or $S^2$. If any is planar, then there is an embedded,
essential planar surface among the fundamental surfaces of $(M,\T)$.
If either $M$ is a knot-manifold or $P$ has all its boundary in a
single component of $\bdy M$, then by the Remark following
Proposition \ref{all-in-one-bdry} and its corollaries, we would
necessarily discover an embedded planar, normal surface among the
fundamental surfaces of $(M,\T)$.

Having made these observations, we consider the first step of a
general algorithm. Given a link-manifold $M$ with $0$--efficient
triangulation $\T$, we first construct the fundamental surfaces of
$(M,\T)$. This gives us a constructible, finite collection,
$\mathcal{F}^0$, of normal surfaces in $M$. We can recognize any of
these surfaces that are planar; and by Theorem
\ref{decide-essential}, we can decide if any of these are essential
in $M$. If $P$ were to exist and if any of the fundamental surfaces
of $(M,\T)$ with projective class in $\C(P)$ were planar, then we
would have found one among the surfaces in the collection
$\mathcal{F}^0$.

Now, we return to our consideration should $P$ exist. From the
preceding paragraphs, we may assume none of the fundamental surfaces
of $(M,\T)$ with projective class in $\C(P)$ are  planar surfaces.
Then it follows that $P$ can be written as a sum
$$P = \sum l_i H^0_i,$$ where $H_i^0$ is fundamental,  $\chi(H^0_i)\le 0$
 and if  $\chi(H^0_i) = 0$, then $H_i^0$ is a torus or Klein bottle.  Set
\[C'_0=max \left\{ \frac{L(\bdy H^0_i)}{-\chi(H^0_i)}\right\},\]
where $\chi(H^0_i)<0$ and is fundamental in $(M,\T)$ with projective
class in $\C(P)$. We have from Corollary \ref{ALE-link2} that
$\lambda^0_{av}<C'_0$, where $\lambda^0_{av}$ is the average length
of the components of $\bdy P$.

Hence, there is a component  of $\bdy P$, having slope (say)
$\alpha^{1}$, with the property that $L(\alpha^1)\le C'_0$, where
$L(\alpha^1)$ is the length of $\alpha^1$ (there is a unique norml
representative of the slope $\alpha^1$). All components of $\bdy P$
in the same component of $\bdy M$ as $\alpha^1$ also have slope
$\alpha^{1}$. Let $M(\alpha^{1})$ denote the triangulated Dehn
filling of $M$ along slope $\alpha^{1}$; $M(\alpha^{1})=
M\cup\bbb{T}(\alpha^{1})$ has a triangulation $\T(\alpha^{1})$ that
is $\T$ on $M$ and a minimal layered-triangulation on the solid
torus $\bbb{T}(\alpha^{1})$. From $P$ we get a planar normal
surface, denoted $P(\alpha^{1})$, in
$(M(\alpha^{1}),\T(\alpha^{1}))$ obtained by ``capping off" the
surface $P$ with copies of the (normal) meridional disk in the solid
torus $\bbb{T}(\alpha^{1})$.

We can write $P(\alpha^1) = \sum \ell^1_j\widehat{H}^1_j$, where the
surfaces $\widehat{H}^1_j$ are fundamental in
$(M(\alpha^{1}),\T(\alpha^{1})) $. But since we are using
triangulated Dehn fillings, it follows from Lemma \ref{re-write}
that each $\widehat{H}^1_j=H^1_j(\alpha^1)$, where $H^1_j$ is a
normal surface in ($M,\T$) and $H^1_j(\alpha^1)$ is $H^1_j$ capped
off with copies of the meridional disk in $\bbb{T}(\alpha^1)$. Thus
$P(\alpha^{1}) = \sum \ell^1_j H^1_j(\alpha^1)$ and we can {\it
re-write} $P$ as $P = \sum \ell^1_j H^1_j$, where $H^1_j$ may not be
fundamental in $M$ but $H^1_j (\alpha^{1})$ is fundamental in
$M(\alpha^1)$. Note that each $H^1_j$ meets the component of $\bdy
M$ containing $\alpha^1$ in the slope $\alpha^1$ and $H^1_j$
projects into $\C(P)$.

For $P(\alpha^{1}) = \sum \ell^1_j H^1_j(\alpha^1)$, if any
$H^1_j(\alpha^1)$ is $S^2, \rppp$ or planar, then  $H^1_j$ is planar
(there are no normal $2$--spheres or projective planes in $M$) and
since $H^1_j$ projects into $\C(P)$ it would be essential. In this
case, we would have an essential, planar surface in $M$ among the
surfaces $\{H^1_j\}$.

Now, we consider the second step of a general algorithm. Having
$\mathcal{F}^0$, the fundamental surfaces in $(M,\T)$, if there are
no essential, planar, fundamental surfaces, we compute
\[C_0=max \left\{ \frac{L(\bdy G^0_i)}{-\chi(G^0_i)}\right\},\]
where $\chi(G^0_i)<0$ and $G^0_i$ is fundamental in $(M,\T)$,
$G^0_i\in \mathcal{F}^0$. Next we find all slopes
$\alpha^1_1,\ldots,\alpha^1_{N_1}$ on $\bdy M$ with
$L(\alpha^1_p)\le C_0$, where $L(\alpha^1_p)$ is the length of a
normal representative of the slope $\alpha^1_p$.

We construct the triangulated Dehn fillings
$(M(\alpha^1_p),\T(\alpha^1_p))$ for all slopes
$\alpha^1_1,\ldots,\\\alpha^1_{N_{1}}$; then in each of these we
compute the fundamental surfaces. From this entire collection of
fundamental surfaces, we select those  in the various
$(M(\alpha^1_p),\T(\alpha^1_p))$ that can be written
$G^1_j(\alpha^1_p)$ for $G^1_j$ a normal surface in $(M,\T)$;
namely, we define $\mathcal{F}^1$ so that $G^1_j \in \mathcal{F}^1$
if and only if there is a fundamental surface $\widehat{G}^1_j$ in
$(M(\alpha^1_p),\T(\alpha^1_p))$, for some $\alpha^1_p$, and
$\widehat{G}^1_j=G^1_j(\alpha^1_p)$, where $G^1_j$ is normal in $M$.
The fundamental surface $\widehat{G}^1_j$ meets the solid torus
$\bbb{T}(\alpha^1_p)$ only in meridional disks and so is some
$G^1_j$ in $M$ ``capped off"; in particular,  $G^1_j$ has a boundary
component with slope $\alpha^1_p$.

 We note
that $C'_0\le C_0$; and therefore, the slope $\alpha^1$ above is a
slope $\alpha^1_p$ for some $p$ and the surfaces $H^1_p$ considered
above are among the surfaces  $G^1_j$.

This gives us a second constructible, finite collection,
$\mathcal{F}^1$, of normal surfaces in $M$. We can recognize any of
these surfaces that are planar and essential in $M$. Again, if $P$
were to exist and if any of the surfaces in $\mathcal{F}^1$ with
projective class in $\C(P)$ were planar, then we would find an
essential planar surface in the collection $\mathcal{F}^1$. Note
that if either $M$ has only two boundary components or $P(\alpha^1)$
has all its boundary in a single component of $\bdy M(\alpha^1)$,
then by the Remark following Proposition \ref{all-in-one-bdry} and
its corollaries, we would necessarily discover an embedded,
essential, planar normal surface in $\mathcal{F}^1$.

Now, suppose $M$ has $n$ boundary components (at this point we may
suppose $n>2$) and $m$ is an integer, $1\le m<n$. We assume we have
determined $m$ collections of slopes on $\bdy M$,
$\{\alpha^1_1,\ldots,\alpha^1_{N_{1}}\},\cdots,
\{\alpha^m_1,\ldots,\alpha^m_{N_{m}}\}$, along with $m+1$
collections of normal surfaces in $(M,\T)$,
$\mathcal{F}^0,\ldots,\mathcal{F}^m$, so that $\mathcal{F}^k$ is the
collection of normal surfaces of $(M,\T)$ and for each $k, 1\le k
\le m$, the normal surface $G^k_i\in \mathcal{F}^k$ if and only if
there is a fundamental surface $\widehat{G}^k_i$ in
$(M(\alpha^1_{j_1},\ldots,\alpha^k_{j_k}),\T(\alpha^1_{j_1},\ldots,\alpha^k_{j_k}))$
for some set of slopes $\alpha^1_{j_1},\ldots,\alpha^k_{j_k}$ with
$\alpha^i_{j_i}\in\{\alpha^i_1,\ldots,\alpha^i_{N_{i}}\}$ and
$\widehat{G}^k_i = G^k_i(\alpha^1_{j_1},\ldots,\alpha^k_{j_k})$.
Furthermore, if there is an embedded, essential, planar, normal
surface in $M$, then either there is an embedded, essential, planar,
normal surface in one of the collections $\mathcal{F}^k$ or there is
an embedded, essential, planar, normal surface $P$ in $M$ that is
least weight in its isotopy class and for some set of slopes
$\{\alpha^1,\ldots,\alpha^m\}$, we can cap off $P$ along these
slopes to get the planar surface $P(\alpha^1,\ldots,\alpha^m)$,
where $\alpha^i\in\{\alpha^i_1,\ldots,\alpha^i_{N_{i}}\}, 1\le i\le
m$.

This gives us a  finite collection,
$\mathcal{F}^0,\mathcal{F}^1,\ldots,\mathcal{F}^m$ of normal
surfaces in $M$. We can recognize among these surfaces any that are
planar and essential in $M$. Again, we observe that if $P$ were to
exist and if any of the surfaces in $\mathcal{F}^i$, $0\le i\le m$,
with projective class in $\C(P)$ were planar, then we would find an
essential planar surface in the collection
$\mathcal{F}^0,\mathcal{F}^1,\ldots,\mathcal{F}^m$. Note that if
either $n = m+1$ ($M(\alpha^1,\ldots,\alpha^m)$ has only one
boundary component) or $P(\alpha^1,\ldots,\alpha^m)$ has all its
boundary in a single component of $\bdy
M(\alpha^1,\ldots,\alpha^m)$, then  we would necessarily discover an
embedded planar, normal surface in
$\mathcal{F}^0,\mathcal{F}^1,\ldots,\mathcal{F}^m$.

We now return to what the situation might be should $P$ exist and we
have not found an embedded, essential, planar surface in the
collection $\mathcal{F}^0,\mathcal{F}^1,\ldots,\mathcal{F}^m$.
Consider $P(\alpha^1,\ldots,\alpha^m)$. We can write
$P(\alpha^1,\ldots,\alpha^m) = \sum \ell^m_k\widehat{H}^m_k$, where
the surfaces $\widehat{H}^m_k$ are fundamental in
$(M(\alpha^1,\ldots,\alpha^m),\T(\alpha^1,\ldots,\alpha^m))$. Again,
since each $\widehat{H}^m_k$ meets the solid tori
$\bbb{T}(\alpha^i)$ in surfaces that sum to the meridional disks, it
follows that each $\widehat{H}^m_k=H^m_k(\alpha^1,\ldots,\alpha^m)$,
where $H^m_k$ is a normal surface in ($M,\T$) and
$H^m_k(\alpha^1,\ldots,\alpha^m)$ is $H^m_k$ capped off with copies
of the meridional disks in the solid tori
$\bbb{T}(\alpha^1),\ldots,\bbb{T}(\alpha^m)$. Thus
$P(\alpha^1,\ldots,\alpha^m) = \sum \ell^m_k
H^m_k(\alpha^1,\ldots,\alpha^m)$ and we can {\it re-write} $P$ as $P
= \sum \ell^m_k H^m_k$, where $H^m_k$ may not be fundamental in $M$
but $H^m_k(\alpha^1,\ldots,\alpha^m)$ is fundamental in
$M(\alpha^1,\ldots,\alpha^m)$. Note that each $H^m_k$ meets the
component of $\bdy M$ containing the slope $\alpha^i$ in that slope
and $H^m_k$ projects into $\C(P)$. In particular, the collection
$\{H^m_k\}$ is a spanning collection for the normal surfaces in $M$
over $\C(P)$, the carrier of $P$.

If any $H^m_k(\alpha^1,\ldots,\alpha^m)$ is $S^2, \rppp$ or planar,
then $H^m_k$ is planar (there are no normal $2$--spheres or
projective planes in $M$) and since $H^m_k$ projects into $\C(P)$ it
would be essential. In this case, we would have an essential, planar
surface in $M$ among the surfaces
$\mathcal{F}^0,\mathcal{F}^1,\ldots,\mathcal{F}^m$.

It follows, from our assumption that there are no essential planar
surfaces in the collections
$\mathcal{F}^0,\mathcal{F}^1,\ldots,\mathcal{F}^m$, that $P$ can be
written as a sum
$$P = \sum l_k H^m_k,$$ where  $\chi(H^m_k)\le 0$
 and if  $\chi(H^m_k) = 0$, then $H_k^m$ is a torus or Klein bottle.

 Set
\[C'_m=max \left\{ \frac{L(\bdy H^m_k(\alpha^1,\ldots,\alpha^m))}{-\chi(H^m_k(\alpha^1,\ldots,\alpha^m))}\right\},\]
where $\chi(H^m_k(\alpha^1,\ldots,\alpha^m))<0$ and is fundamental
with projective class in $\C(P(\alpha^1,\\\ldots,\alpha^m))$. We
have from Corollary \ref{ALE-link2} that $\lambda^m_{av}\le C'_m$,
where $\lambda^m_{av}$ is the average length of the components of
$\bdy P(\alpha^1,\ldots,\alpha^m)$.

Hence, there is a component of $\bdy P(\alpha^1,\ldots,\alpha^m)$
having slope (say) $\alpha^{m+1}$ with the property that
$L(\alpha^{m+1})\le C'_m$, where $L(\alpha^{m+1})$ is the length of
$\alpha^{m+1}$. All components of $\bdy P(\alpha^1,\ldots,\alpha^m)$
in the same component of $\bdy M(\alpha^1,\ldots,\alpha^m)$ as
$\alpha^{m+1}$ also have slope $\alpha^{m+1}$. Let
$M(\alpha^1,\ldots,\alpha^m,\alpha^{m+1})$ denote the triangulated
Dehn filling of $M(\alpha^1,\ldots,\alpha^m)$ along slope
$\alpha^{m+1}$; $M(\alpha^1,\ldots,\alpha^m,\alpha^{m+1})=
M(\alpha^1,\ldots,\alpha^m)\cup\bbb{T}(\alpha^{m+1})$ has a
triangulation $\T(\alpha^1,\ldots,\alpha^m,\alpha^{m+1})$ that is
$\T$ on $M$ and a minimal layered-triangulation on each
$\bbb{T}(\alpha^{i}),i=1,\ldots,m+1$. From $P$ we get a planar
normal surface, denoted $P(\alpha^1,\ldots,\alpha^m,\alpha^{m+1})$,
in
$(M(\alpha^1,\ldots,\alpha^m,\alpha^{m+1}),\T(\alpha^1,\ldots,\alpha^m,\\\alpha^{m+1}))$
obtained by ``capping off" the surface $P$ with copies of the
(normal) meridional disks in the solid tori $\bbb{T}(\alpha^{i})$.

We can write $P(\alpha^1,\ldots,\alpha^m,\alpha^{m+1}) = \sum
\ell^{m+1}_r\widehat{H}^{m+1}_r$, where the surfaces
$\widehat{H}^{m+1}_r$ are fundamental in
$(M(\alpha^1,\ldots,\alpha^m,\alpha^{m+1}),\T(\alpha^1,\ldots,\alpha^m,\alpha^{m+1}))
$. But just as above, each
$\widehat{H}^{m+1}_r=H^{m+1}_r(\alpha^1,\ldots,\alpha^m,\alpha^{m+1})$.
 Thus we {\it re-write} $P$ as
$P = \sum \ell^{m+1}_rH^{m+1}_r$, where $H^{m+1}_r$ may not be
fundamental in $M$ but $H^{m+1}_r
(\alpha^1,\ldots,\alpha^m,\alpha^{m+1})$ is fundamental in
$M(\alpha^1,\ldots,\alpha^m,\alpha^{m+1})$. Note that each
$H^{m+1}_r$ meets the component of $\bdy M$ containing the slope
$\alpha^{m+1}$ in the slope $\alpha^{m+1}$ and $H^{m+1}_r$ projects
into $\C(P)$.

If any $H^{m+1}_r(\alpha^1,\ldots,\alpha^m,\alpha^{m+1})$ is $S^2,
\rppp$ or planar, then $H^{m+1}_r$ is planar and as earlier, we
would have an essential, planar surface in $M$ among the surfaces
$H^{m+1}_r$.

Finally, this brings us to the completion of the induction step,
going from $m$ to $m+1$ in our construction of the family
$\mathcal{F}^{m+1}$, having the families
$\mathcal{F}^0,\ldots,\mathcal{F}^m$.

Compute
\[C_m=max \left\{ \frac{L(\bdy G^m_r(\alpha^1_{j_1},\ldots,\alpha^m_{j_m}))}{-\chi(G^m_r(\alpha^1_{j_1},\ldots,\alpha^m_{j_m}))}\right\},\]
where $\chi(G^m_r(\alpha^1_{j_1},\ldots,\alpha^m_{j_m}))<0$ and
$G^m_r(\alpha^1_{j_1},\ldots,\alpha^m_{j_m})$ is fundamental in
$(M(\alpha^1_{j_1},\\\ldots,\alpha^m_{j_m}),\T(\alpha^1_{j_1},\ldots,\alpha^m_{j_m}))$
for some set of slopes $\alpha^1_{j_1},\ldots,\alpha^m_{j_m}$ with
$\alpha^i_{j_i}\in\{\alpha^i_1,\ldots,\\\alpha^i_{N_{i}}\}$. Next
find all slopes, say $\alpha^{m+1}_1,\ldots,\alpha^{m+1}_{N_{m+1}}$,
on $\bdy M(\alpha^1_{j_1},\ldots,\alpha^m_{j_m})$, for each set of
slopes $\alpha^1_{j_1},\ldots,\alpha^m_{j_m}$ with
$\alpha^i_{j_i}\in\{\alpha^i_1,\ldots,\alpha^i_{N_{i}}\}$ and having
 $L(\alpha^{m+1}_j)\le C_m$, where $L(\alpha^{m+1}_j)$ is the
length of a normal representative of the slope $\alpha^{m+1}_j$.

Construct the triangulated Dehn fillings
$(M(\alpha^1_{j_1},\ldots,\alpha^m_{j_m},\alpha^{m+1}_{j_{m+1}}),
\T(\alpha^1_{j_1},\ldots,\alpha^m_{j_m},\\\alpha^{m+1}_{j_{m+1}}))$
for some set of slopes
$\alpha^1_{j_1},\ldots,\alpha^{m+1}_{j_{m+1}}$ with
$\alpha^i_{j_i}\in\{\alpha^i_1,\ldots,\alpha^i_{N_{i}}\}$, all $i,
1\le i\le m+1$. Having these triangulated Dehn fillings, we find for
each their fundamental surfaces. We define $\mathcal{F}^{m+1}$ so
that $G^{m+1}_s \in \mathcal{F}^{m+1}$ if and only if there is a
fundamental surface $\widehat{G}^{m+1}_s$ in one of our Dehn
fillings
$(M(\alpha^1_{j_1},\ldots,\alpha^m_{j_m},\alpha^{m+1}_{j_{m+1}}),
\T(\alpha^1_{j_1},\ldots,\alpha^m_{j_m},\\\alpha^{m+1}_{j_{m+1}}))$
and
$\widehat{G}^{m+1}_s=G^{m+1}_s(\alpha^1_{j_1},\ldots,\alpha^m_{j_m},\alpha^{m+1}_{j_{m+1}})$,
where $G^{m+1}_s$ is normal in $M$.

We have $C'_m\le C_m$; and therefore, the slope $\alpha^{m+1}$ above
is a slope $\alpha^{m+1}_j$ for some $j$ and the surfaces
$\{H^{m+1}_r\}$ considered above are among the surfaces
$\{G^{m+1}_s\}$.

This defines the finite family $\mathcal{F}^{m+1}$ of normal
surfaces in $M$. We can recognize any of these surfaces that are
planar and essential in $M$. Again, if $P$ were to exist and if any
of the surfaces in $\mathcal{F}^{m+1}$ with projective class in
$\C(P)$ were planar, then we would find an essential planar surface
in the collection $\mathcal{F}^{m+1}$.

If there is a properly embedded, essential, planar surface  in $M$,
then we will find one in one of the collections
$\mathcal{F}^0,\ldots\mathcal{F}^{n-2}$ or, by Corollary
\ref{knot-planar}, we must find one in the collection
$\mathcal{F}^{n-1}$, since each
$M(\alpha^1_{j_1},\ldots,\alpha^{n-1}_{j_{n-1}})$ is a
knot-manifold. If we do not find an essential, planar surface in any
of the collections $\mathcal{F}^i, 0\le i\le (n-1)$, then there is
no properly embedded, essential, planar surface in $M$.\end{proof}

Just a final remark, which was alluded to at various points
throughout the proof. We commented earlier that in writing a planar
normal surface $P$ as a sum of normal surfaces, $F_1,\ldots,F_k$,
then it is quite possible that none of the $F_i$ are planar.
However, what the argument above accomplishes is that we keep
re-writing $P$ as a geometric sum of normal surfaces that agree with
the slope of the curves in $\bdy P$ on more and more of the boundary
tori in $M$, while each new summand still projects into the carrier
of $P$ in $\P(M,\T)$. Finally, at least if we finally get to a knot
manifold, we can re-write $P$ as a sum where for $b_i$ the number of
boundary components of $F_i$ and $b$ the number of boundary
components of $P$, we have $b=\sum n_ib_i$, ``the boundary of the
sum is the sum of the boundaries". In this situation, if $P$ is
planar, then one of the $F_i$ must also be planar. Such a planar
$F_i$ is essential because it projects into $\C(P)$.

\section{Finding planar surfaces with boundary conditions}
As mentioned earlier, our original interest in studying planar
surfaces in knot- and link-manifolds was an interest in singular
planar surfaces and, specifically, their application to problems
like the Word Problem and Conjugacy Problem. These problems have
boundary conditions. Hence, in hope of better understanding the
problem for singular surfaces, our considerations evolved to
analogous questions about embedded, planar surfaces. There are four
such questions, with varying types of boundary conditions.

Recall that if $\mu$ and $\lambda$ are slopes in a torus, we use
$\langle\mu,\lambda\rangle$ to denote their distance; and if $\mu$
is a given slope and $\lambda$ is a slope with
$\langle\mu,\lambda\rangle = 1$, we call $\lambda$ a {\it longitude
with respect to $\mu$}. Sometimes when we want to distinguish a
particular slope in a component of the boundary of a given knot- or
link-manifold, we call it a meridian. In particular, this is the
case when we are making statements analogous to the singular
problems related to the Word or Conjugacy Problems.

Suppose we are given a link-manifold $M$ and a component $B$ of
$\bdy M$. The various boundary conditions we consider are:

\begin{enumerate} \item[1.] {\bf (Condition $\boldsymbol{B}$)} Can it be decided if
there is a properly embedded, essential, planar surface in $M$ with
boundary meeting $B$?\item[2.]{\bf (Condition
$\boldsymbol{B}_\gamma$)} Given a slope $\gamma$ in $B$, can it be
decided if there is a properly embedded, essential, planar surface
in $M$ meeting $B$ in the slope $\gamma$?\item[3.]{\bf (Condition
$\boldsymbol{E}_\gamma$)} Given a slope $\gamma$ in $B$, can it be
decided if there is a properly embedded, essential {\it
punctured-disk} in $M$ with boundary slope $\gamma$ and punctures in
$\bdy M\setminus B$?\item[4.]{\bf(Condition $\boldsymbol{E}$)} Given
a meridional slope $\mu$ in $B$, can it be decided if there is a
properly embedded, essential {\it punctured-disk} with boundary a
longitude in $B$ and punctures in $\bdy M\setminus
B$?\end{enumerate}

We show there are algorithms for  Condition $E_\gamma$ and Condition
$E$; however, we are not able to answer the question completely for
either Condition $B$ or Condition $B_\gamma$, except in the case of
knot-manifolds. Condition $B$ seems, on first glance, much like the
question of the previous section. However, a closer look shows that
in the previous section, we might have an essential planar surface
with its boundary in $B$ but the planar surface produced from it,
which had to be in our collection of constructible surfaces to
check, did not have its boundary in $B$. There are examples that
clearly expose the deficiency of the method in the previous section.
Namely, the planar surface $P$, with boundary meeting $B$, is a
geometric sum of a collection of planar surfaces none of which have
boundary in $B$ and a higher genus surface with boundary in $B$.
Under geometric sum, the handles for the higher genus surface are
replaced by copies of the planar surfaces, reducing genus but adding
boundary to get $P$. Such a sum could be quite complicated and it is
not clear, if this were the case, how to find $P$.

However, without modification of the proof of Theorem
\ref{find-planar}, we have the following theorem.
\begin{thm}\label{modify-B} Given a link-manifold $M$ and a component $B$ of $\bdy
M$, there is an algorithm to decide if there is an embedded,
essential, planar surface either with boundary meeting $B$ or with
no boundary meeting $B$; or neither of these possibilities occur. If
a planar surface exists, the algorithm will construct one.\end{thm}

This result, of course, could have been stated earlier and then we
would have Theorem \ref{find-planar} as its corollary.

If we are also given a slope $\gamma$ on $B$, as in Condition
$B_\gamma$, then, again, without modification of the proof of
Theorem \ref{find-planar} and using the notion given below on normal
surfaces with boundary conditions, the conclusion of Theorem
\ref{modify-B} can be replaced with:

 {\it there is an algorithm to
decide if there is an embedded, essential planar surface either with
slope $\gamma$ on $B$ or with no boundary on $B$; or neither of
these possibilities occur. If a planar surface exists, the algorithm
will construct one.}

\subsection{Boundary conditions on normal surfaces.} The remaining problems
we investigate in this section search for normal surfaces with a
given boundary slope; this is the prototypical example of a boundary
condition on a normal surface.

In what follows we assume the reader is familiar with normal curves
in triangulated $2$--manifolds.

\begin{figure}
\psfrag{m}{\footnotesize{$a_m$}}
\psfrag{l}{\footnotesize{$a_l$}}\psfrag{k}{\footnotesize{$a_k$}}\psfrag{a}{\large$\boldsymbol{\sigma}$}
\psfrag{b}{\large$\boldsymbol{\beta}$}
\psfrag{i}{$\delta_i$}\psfrag{j}{$\delta_j$}\psfrag{1}{$\vec{\delta}_i\rightarrow
\vec{a}_l+\vec{a}_m$}\psfrag{2}{$\vec{\delta}_j\rightarrow\vec{a}_k$}{\epsfxsize
= 1.5 in \centerline{\epsfbox{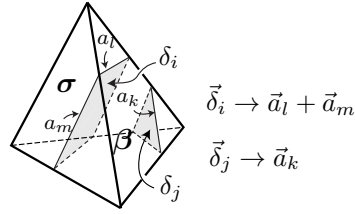}} } \caption{Normal
boundary operator in the case that the faces $\sigma$ and $\beta$
are in $\bdy M$.} \label{f-bound-oper}
\end{figure}

Suppose $M$ is a $3$--manifold and $\T$ is a triangulation of $M$.
Select some order for the normal triangle and quadrilateral types in
the tetrahedra of $\T$, which then determines coordinates in
$\bbb{R}^{7t}$; similarly, select some order for the normal arc
types in $\bdy M$, which then determines coordinates in $\bbb{R}^m$,
where $m$ is the number of arc types in the triangulation induced on
$\bdy M$. Now, if $\delta_i$ is a normal triangle or quad type in
$M$, we associate with $\delta_i$ the unit vector
$\vec{\delta}_i=(x_1,\ldots,x_q,\ldots,x_{7t})$, where $x_q=0,q\ne
i, x_i=1$; similarly, if $a_k$ is a normal arc type in $\bdy M$, we
associate the unit vector $\vec{a}_k = (z_1,\dots,z_p,\ldots,z_m)$,
where $z_p=0,p\ne k, z_k=1$. There is a {\it normal boundary
operator}, denoted $\bdy_{\T}$, which is a linear map from
$\bbb{R}^{7t}$ to $\bbb{R}^m$ defined as
$$\bdy_{T}(\vec{\delta}_i)=\sum \epsilon_{i,j}\vec{a}_j,$$ where
$\epsilon_{i,j} = 1$ if the arc type $a_j$ in $\bdy M$ is in the
normal disk type $\delta_i$ and $\epsilon_{i,j} = 0$, otherwise. See
Figure \ref{f-bound-oper}.

Now, suppose $M$ is a link-manifold, $B$ is a boundary component and
$\T$ is a triangulation of $M$ that induces a one vertex
triangulation on $B$. From \cite{jac-sedg-dehn}, we have that the
slopes on $B$ are in one-one correspondence with normal isotopy
classes of nontrivial simple closed curves in $B$. In the above
scheme, it follows that if $\gamma$ is a slope, then $\gamma$
corresponds to a unique integer lattice point in the parametrization
space, $\bbb{R}^m$, of normal curves in $\bdy M$ and thus determines
a unique linear subspace, say $R(\gamma)$. The inverse of
$R(\gamma)$ under the boundary operator $\bdy_{\T}$ is a vector
subspace of $\bbb{R}^{7t}$ and meets the normal solution space,
$\S(M,\T)$, for normal surfaces in $M$ in a sub-cone. This sub-cone,
contains the parametrization of those normal surfaces in $M$ that
either meet $B$ only in the slope $\gamma$ or do not meet $B$ at all
(which includes the closed normal surfaces). We denote this sub-cone
$\S_\gamma(M,\T)$. It satisfies the algebra we developed for the
normal solution space, such as having a projective solution space,
which will be a subcell of the projective space for $M$, having
vertex and fundamental solutions associate with it (however, these
are not necessarily vertex and fundamental solutions of $\S(M,\T)$),
and so on.

\begin{figure}
\psfrag{1}{$a_1$}
\psfrag{2}{$a_2$}\psfrag{3}{$a_3$}\psfrag{4}{$a_4$}\psfrag{5}{$a_5$}\psfrag{6}{$a_6$}
\psfrag{B}{$B$}\psfrag{a}{$z_1=z_4$}\psfrag{b}{$z_2=z_5$}\psfrag{c}{$z_3=z_6$}
\psfrag{x}{\Large{$\frac{z_1}{z_3} = \frac{2}{5}$}}{\epsfxsize = 3.5
in \centerline{\epsfbox{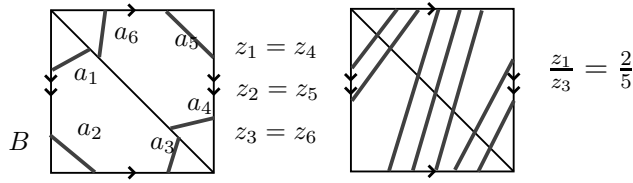}} } \caption{Normal arcs
in the one-vertex triangulation of the torus, the equations for
normal curves, and an example of the slope corresponding to the
ratio $z_1/z_3 = 2/5$.} \label{F-BOUND-planar}
\end{figure}

While the above description of $\S_\gamma(M\T)$ provides a nice
theoretical setting, there is a direct method of determining
$\S_\gamma(M,\T)$ using a system of homogeneous linear equations. To
see this recall from \cite{jac-sedg-dehn} that since $\T$ induces a
one-vertex triangulation on $B$, the normal isotopy classes of
normal curves in $B$ can be parameterized by just three normal arc
types in one of the two triangles in the induced triangulation on
$B$ ($m$ above can be taken to be $3$). See Figure
\ref{F-BOUND-planar}. Choose notation so that the three arc types
are $a_1, a_2$ and $a_3$. If we consider only essential curves in
$B$ and prescribe a slope $\gamma$, then we have that one of the
coordinates must be zero (otherwise, we would allow trivial curves)
and the other two must be in constant ratio. In particular, if $z_k,
k=1,2, 3$, is the number of arcs of type $a_k$ in $\gamma$ and, say,
$a_2 = 0$, then $z_1, z_3$ are in a constant ratio; i.e., there are
relatively prime integers $r,s$ (possibly $r=1=s$) so that $z_1/z_3
= r/s$. Denote the normal triangles that meet $B$ in the arc type
$a_k$ by $\delta_{i_k}$ and the normal quadrilaterals that meet $B$
in the arc type $a_k$ by $\delta^Q_{j_k}$. Then if we add to the
matching equations for $(M,\T)$ the equations
$$x_{i_2} + y_{j_2}=0$$ and $$ s(x_{i_1}+y_{j_1}) = r(x_{i_3}+y_{j_3}),$$ we recover
$\S_\gamma(M,\T)$ as the solution cone in the positive orthant of
$\bbb{R}^{7t}$ for this system of equations.

\subsection{Boundary conditions and average length estimate.}
Suppose $M$ is a link-manifold, $B$ is a component of $\bdy M$, and
$\gamma$ is a slope on $B$. Now, suppose $F$ is an embedded normal
surface and  $b_\gamma$ is the number of components of $\bdy F$
having slope $\gamma$. Let $b_{-\gamma}$ denote the number of
components of $\bdy F$ that do not have slope $\gamma$; if $b$ is
the number of components of $\bdy F$, $b = b_\gamma + b_{-\gamma}$.
If $b_{-\gamma}\ne 0$, we set $\lambda^{-\gamma}_{av}= [(L(\bdy
F)-b_\gamma L(\gamma)]/b_{-\gamma}$ and call
$\lambda^{-\gamma}_{av}$ the {\it average length of $\bdy F$ away
from $\gamma$}. If $F$ does not meet $B$ in slope $\gamma$, then
this is just the average length of $\bdy F$. We have the following
modification of Proposition \ref{ALE-link1}

\begin{prop}\label{ALE-link-bound} Suppose $M$ is a
link-manifold with no embedded annuli having boundary in distinct
components of $\bdy M$, $B$ is a component of $\bdy M$, and $\gamma$
is a slope on $B$. Furthermore, suppose $\T$ is a $0$--efficient
triangulation of $M$. Then there is a constant $C=C(M,\T)$,
depending only on $M$ and $\T$ so that if $F$ is an embedded normal
surface in $M$ with no trivial boundary curves and
$\lambda^{-\gamma}_{av}$ is the average length of the components of
$\bdy F$ away from $\gamma$, then
$$\lambda^{-\gamma}_{av}\le C(2g+2+b_\gamma),$$ where $g$ is the genus
of $F$ and $b_\gamma$ bounds the number of components of $\bdy F$
having slope $\gamma$.\end{prop}

\begin{proof} We use the same notation and proceed just as in the
proof of Proposition \ref{ALE-link1}. We define $$C =
max\left\{\frac{L(\bdy F_1)}{-\chi(F_1)},\ldots,\frac{L(\bdy
F_I)}{-\chi(F_I)}, \frac {L(\bdy A_1)}{\abs{A_1}},\ldots,\frac
{L(\bdy A_K)}{\abs{A_K}}\right\}.$$

We have $b_\gamma$ bounding the number of components of $\bdy F$
having slope $\gamma$; $b_{-\gamma}$ the number of components of
$\bdy F$ not having slope $\gamma$; $b_A = \sum n_k\abs{A_k}$, the
number of boundary components of $\sum n_kA_k$; and $b$ the number
of boundary components of $F$. It follows that
$$b_{-\gamma}\lambda^{-\gamma}_{av}\le b_{-\gamma}\lambda^{-\gamma}_{av}+
(b-b_{-\gamma})L(\gamma) = L(\bdy F) = \sum l_iL(\bdy F_i) + \sum
n_kL(\bdy A_k)\le$$ $$\le C\left(\sum l_i(-\chi(F_i))+\sum
n_k\abs{\bdy A_k}\right)=C(-\chi(F') + b_A).$$  However, $\chi(F) =
\chi(F')$; hence,  we have $b_{-\gamma}\lambda^{-\gamma}_{av}\le
C(2g-2+b+b_A)$. It follows that $$\lambda^{-\gamma}_{av}\le
C\left(\frac{2g-2}{b_{-\gamma}}
+\frac{b_\gamma}{b_{-\gamma}}+1+\frac{b_A}{b_{-\gamma}}\right) \le
C\left(2g +\frac{b}{b_{-\gamma}}+1\right)\le C(2g+
2+b_\gamma).$$\end{proof}

\subsection{Conditions $E_\gamma$ and  $E$.} We first construct an
algorithm for Condition $E_\gamma$, using induction on the number of
boundary components of the given link-manifold. Afterwards, given a
link-manifold, we are able to derive an algorithm for Condition $E$
by using the algorithm for Condition $E_\gamma$ in a related
link-manifold, which is constructed from the given link-manifold.
For our induction step in the algorithm for $E_\gamma$, we use Dehn
filling to reduce the number of boundary components and invoke the
induction hypothesis. However, we essentially go the opposite
direction in deriving the algorithm for Condition $E$; in this case,
we construct a new link-manifold by removing a knot from the given
link-manifold, a process which we might think of as ``Dehn
drilling".

\begin{figure}[htbp]

           \psfrag{i}{\large{$\eta = D\times S^1$}}
\psfrag{a}{$A$} \psfrag{1}{$B_1$} \psfrag{2}{$B_2$}
\psfrag{3}{$B_3$}
\psfrag{P}{\large{$P$}}\psfrag{b}{$bd(P)$}\psfrag{p}{punctures}

        \vspace{0 in}
        \begin{center}
\epsfxsize = 3.5 in \epsfbox{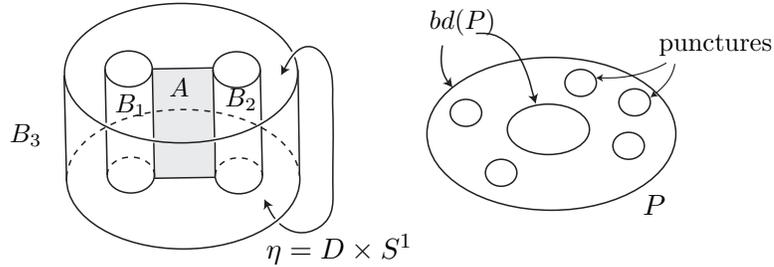}
 \caption{A pair-of-pants crossed with the circle, $\eta$, which is
homeomorphic with the neighborhood $\eta(A)$ of $B'\cup A\cup B''$
where $A$ is an essential annulus in a link-manifold having boundary
in distinct boundary components $B'(B_1)$ and $B''(B_2)$. On the
right is a punctured-annulus $P$ with boundary $bd(P)$.}
\label{f-pair-of-pants-times}
\end{center}
\end{figure}

Suppose $D$ is a planar surface with three boundary components (a
pair-of-pants). We wish to understand the embedded, essential planar
surfaces in $\eta = D\times S^1$. Suppose $b_1, b_2$ and $b_3$ are
the boundary components of $D$ and $B_i=b_i\times S^1, i = 1, 2, 3$,
denotes the corresponding components of $\bdy\eta$. For $t\in S^1$,
we call the slope determined by $t\times b_i$ meridional on $B_i$
and denote it by $\mu_i$;  for any point $x\in b_i$, we call the
slope determined by $x\times S^1$ vertical and denote it by
$\alpha_i$. The product structure on $\eta$ makes it a Seifert fiber
space. See Figure \ref{f-pair-of-pants-times}.

Just as we distinguished a special boundary component of a planar
surface in our definition of a punctured-disk, if we have a planar
surface with at least two boundary components and we distinguish two
distinct boundary components, we call the planar surface a {\it
punctured-annulus} and refer to all boundary components distinct
from the two we have distinguished as the punctures. Again, see
Figure \ref{f-pair-of-pants-times}. If $P$ is a punctured-annulus,
$C'$ and $C''$ the two boundary components we have distinguished,
then we set $bd(P) = C'\cup C''$ and call $bd(P)$ the boundary of
$P$.

For each essential arc $a$ in $D$, we get an embedded, essential,
vertical annulus $a\times S^1$ in $\eta$. There are three that have
their boundaries in distinct components of $\bdy \eta$; we denote
these $A_{1,2}, A_{1,3}$ and $A_{2,3}$, where $A_{i,j}$ is a
vertical annulus having one boundary, $\alpha_i$ in $B_i$ and the
other $\alpha_j$ in $B_j$.

Since there is a symmetry between the boundary components of $\eta$,
we shall distinguish one, say $B_2$, and describe a family of
embedded, essential planar surfaces in $\eta$. Each of the surfaces
we describe is a punctured-annulus having one boundary in $B_1$, the
other in $B_3$, and all its punctures in $B_2$. We get an analogous
family by distinguishing either $B_1$ or $B_3$. Notice that a Dehn
filling of $\eta$ along the slope $\mu_2$ gives the manifold
$\eta(\mu_2) = (S^1\times S^1)\times [0,1]$, a torus cross the
interval, with $B_ 1= S^1\times S^1\times \{1\}$ and $B_3 =
S^1\times S^1\times \{0\}$. For every slope $\gamma_1$ in $B_1$
there is an annulus $\gamma_1\times [0,1]$ in $\eta(\mu_3)$ meeting
$B_1$ in $\gamma_1$ and $B_3$ in a corresponding slope, say
$\gamma_3$; and from this annulus, we obtain a punctured-annulus in
$\eta$, having all its punctures in $B_2$ in the slope $\mu_2$. We
denote this punctured-annulus $P^\gamma_{1,3}$. Now, having the
various vertical annuli, $A_{i,j}$, described above, we can Dehn
twist the punctured-annulus $P^\gamma_{1,3}$ about $A_{i,j}$. Dehn
twisting about $A_{1,3}$ gives the same punctured-annuli that we
would have gotten from our construction had we first twisted the
slope $\gamma_1$ about the vertical slope $\alpha_1$ on $B_1$.
However, Dehn twisting about $A_{1,2}$ or $A_{2,3}$ gives  new
families of punctured-annuli. By twisting $P^\gamma_{1,3}$ about
$A_{1,2}$ we get a family of punctured-annuli that meet  $B_1$ in
the family of slopes obtained from  $\gamma_1$ by Dehn twisting
about the vertical slope $\alpha_1$, using the same number of twist
as we have about $A_{1,2}$; they meet $B_3$ in the slope $\gamma_3$;
and they meet $B_2$ in a number of punctures, each having a slope
obtained from $\mu_2$ by Dehn twisting about $\alpha_2$, using the
same number of twist as we have about $A_{1,2}$. We get an analogous
family if we twist about $A_{2,3}$ except that we do not change the
slope on $B_1$ in these cases.

\begin{lem}\label{planar-pair-pants} The above examples describe all
possible embedded, essential planar surfaces in $\eta = D\times
S^1$. \end{lem}
\begin{proof} Suppose $P$ is an embedded, essential planar surface with boundary
in only one component of $\bdy\eta$, say $B_1$. Then by Dehn filling
along the meridian slopes, $\mu_2$ and $\mu_3$, we get a solid
torus, $\bbb{T} = \eta(\mu_2,\mu_3)$. The only planar surfaces in
$\bbb{T}$ meeting $B_1 = \bdy \bbb{T}$ in essential curves are the
meridional disk and an annulus parallel into $B_1$. But $P$ having
boundary only in $B_1$ leaves only the possibility that $P$ is a
vertical annulus; furthermore, to be essential in $\eta$, $P$ must
separate the boundary components, $B_2$ and $B_3$.

If $P$ meets precisely two of the boundary components, say $B_1$ and
$B_2$, then we can Dehn fill $\eta$ along $\mu_3$. The Dehn filling
$\eta(\mu_3)$ is an annulus cross the circle (torus cross an
interval) and all embedded, essential planar surfaces with their
boundaries in distinct components of the boundary are annuli. For
such a surface not to meet the third boundary component of $\eta$,
it must be vertical in $\eta$.

The only remaining possibility is that $P$ meets all three boundary
components of $\eta$. In this case, we note that $P$ can not meet a
boundary component in its vertical slope. If $P$ did meet, say $B_3$
in $\alpha_3$, then we can Dehn fill $\eta$ along $\alpha_3$,
getting the essential planar surface $P(\alpha_3)$ obtained by
capping off $P$ in $\eta(\alpha_3)$, which is the connected sum of
two solid tori. But since $P(\alpha_3)$ has boundary in each
boundary component of $\eta(\alpha_3)$, $P$ could not be essential
in $\eta$.

It follows that after Dehn filling $\eta$ along all the boundary
slopes of $P$, the Seifert fiber structure on $\eta$ extends to a
Seifert fiber structure on the Dehn filled manifold; however,
capping off $P$ gives an embedded, horizontal $2$--sphere. It
follows that one of the slopes of $\bdy P$ on the boundary of $\eta$
is meridional (possibly after a Dehn twist about a vertical
annulus). Thus $P$ is a punctured-annulus with its punctures on this
boundary component of $\eta$ that $P$ meets in meridional slopes.
Reversing the Dehn twist, gives us one of the punctured-annuli
described above.\end{proof}

Splitting a manifold along a properly embedded surface is a standard
notion; however, we want a special form of this when splitting a
link-manifold along a properly embedded annulus having its boundary
in distinct components of the boundary of the link-manifold. Suppose
$M$ is a link-manifold with $B'$ and $B''$ distinct boundary
components of $\bdy M$; and suppose $A$ is an embedded, essential
annulus with one boundary in $B'$ and one boundary in $B''$. Denote
a small regular neighborhood of $B'\cup A\cup B''$ by $\eta(A)$.
Then $\eta(A)$ is a disk with two punctures (a pair of pants)
crossed with the circle; see Figure \ref{f-pair-of-pants-times}. If
we denote the frontier of $\eta(A)$ by $B_A$, then the boundary
components of $\eta(A)$ are $B',B''$, and $B_A$. Let $M_A =
M\setminus\open{\eta}(A)$, then $M_A$ is a link-manifold with a
boundary component $B_A$ and one fewer boundary component than the
link-manifold $M$. We say $M_A$ is obtained from $M$ by {\it
splitting $M$ along $A$}. Finally, since $\eta(A)$ is a
pair-of-pants cross $S^1$,  a slope on any one of the three boundary
components can be uniquely associate with a ''parallel" slope on the
other two boundary components; we shall use the convention in this
situation of $\gamma'$ being the slope on $B'$, $\gamma''$ on $B''$
and $\gamma_A$ on $B_A$ (if $B'=B$, then we use $\gamma$ for
$\gamma'$). We refer to the meridional slopes designated above for
$\eta$ as meridional slopes on $B'$ ($B$), $B''$ or $B_A$ in
$\eta(A)$ and write $\mu_A$ for $B_A$, and so on. In Lemma
\ref{planar-pair-pants}, we characterized embedded, essential planar
surfaces in $\eta(A)$. Now, we  provide the relationship between
embedded, essential planar surfaces in $M$ and those in $M_A$.

\begin{prop}\label{planar-induct-step} Suppose $M$ is a
link-manifold, $B$ is a component of $\bdy M$ and $A$ is an
embedded, essential annulus having its boundary in distinct
components $B'$ and $B''$ of $\bdy M$. Let $M_A$ be the
link-manifold obtained by splitting $M$ along $A$. We have the
following:\begin{enumerate}\item For $B'\ne B\ne
B''$.\begin{enumerate}\item There is an embedded, essential planar
surface in $M$ with boundary meeting $B$ if and only if there is one
in $M_A$ with boundary meeting $B$.
\item  There is an embedded punctured-disk in $M$ with its
boundary meeting $B$ in slope $\gamma$ if and only if there is one
in $M_A$ with its boundary meeting $B$ in slope
$\gamma$.\end{enumerate}\item For $B'= B \ne
B''$.\begin{enumerate}\item Either the only embedded, essential
planar surfaces meeting $B$ are annuli meeting in slope $\alpha$ and
there are no embedded, essential planar surfaces in $M_A$ meeting
$B_A$ or there is an embedded, essential planar surface in $M$ with
boundary meeting $B$ if and only if there is one in $M_A$ with
boundary meeting $B_A$.\item There are embedded punctured-disks in
$M$ with boundary meeting $B$ in every slope if and only if there is
an embedded punctured-disk in $M_A$ with boundary meeting $B_A$ in a
slope obtained from $\mu_A$ by Dehn twisting about $\alpha_A$ (i.e.,
a slope having geometric intersection one with $\alpha_A$).\item If
the previous situation does not hold, either the only embedded
punc-tured-disk in $M$ meeting $B$ is $A$ and there are no embedded
punctured-disks in $M_A$ meeting $B_A$ or there is an embedded
punctured-disk in $M$ with boundary meeting $B$ in slope $\gamma$ if
and only if there is an embedded punctured-disk in $M_A$ with
boundary meeting $B_A$ in a slope obtained from $\gamma_A$ by Dehn
twisting about $\alpha_A$.\end{enumerate}
\end{enumerate}
\end{prop}

\begin{proof} Proof for 1(a)and (b). Suppose $P$ is an
embedded, essential planar surface in $M$ and its boundary meets
$B$. We isotope $P$ so that it meets $B_A$ minimally and
transversely; this does not affect any components of $\bdy P$ in
$B$. Hence, if $P$ meets $\eta(A)$ at all, it must meet $\eta(A)$ in
embedded, essential planar surfaces in $\eta(A)$. It follows from
Lemma \ref{planar-pair-pants} that removing such pieces from $P$
leaves a connected embedded, essential planar surface in $M_A$
meeting $B$ exactly as $P$ did.

Conversely, suppose $P_A$ is an embedded, essential planar surface
in $M_A$ and it meets $B$. If it does not meet $B_A$, it is also an
embedded, essential planar surface in $M$. If it does meet $B_A$,
then we consider the slope of its boundary in $B_A$. By Lemma
\ref{planar-pair-pants}, no matter the slope there is an embedded,
essential punctured annulus in $\eta_A$ having just one of its
boundary components in $B_A$ and having this slope. If $P_A$ meets
$B_A$ in $m$ components, we add $m$ copies of such a
punctured-annulus in $\eta_A$, arriving at the desired planar
surface for $M$. Notice that this also proves  part 1(b) as well,
since none of this had any affect on the meets with $B$.

Proof of 2(a). Here we might have an essential, planar surface in
$M$ meeting $B$ but there are none in $M_A$ meeting $B_A$; that is,
the planar surface in $M$ is contained in $\eta(A)$. However, from
Lemma \ref{planar-pair-pants}, we have that such a surface is then
one of the vertical annuli in $\eta(A)$. Otherwise, there is an
embedded essential surface in $M_A$ meeting $B_A$ if and only if
there is one in $M$ meeting $B$.

Proof of 2(b) and (c). Suppose we have an embedded punctured-disk in
$M$ with its boundary in $B$ having slope $\gamma$. Then by an
isotopy, we may make it meet $B_A$ transversely and minimally. Thus
we have that a component of its intersection with $\eta(A)$ is a
punctured-disk in $\eta(A)$ with its boundary having slope $\gamma$
in $B$ (if it has more than one component in $\eta(A)$, then each
must be a vertical annulus and have at least one boundary component
in $B_A$; in which case $\gamma = \alpha$). We use the
characterization of planar surfaces in $\eta(A)$ given in Lemma
\ref{planar-pair-pants} to determine the possibilities.

If our punctured-disk meets $\eta(A)$ only in vertical annuli, then
$\gamma = \alpha$ and we have the possibility that there are no
punctured-disk in $M_A$ having boundary in $B_A$ or there is a
punctured-disk in $M_A$ having its boundary of slope $\alpha_A$ in
$B_A$. This satisfies the conclusion to part 2(c). Notice that if
there are no punctured-disk in $M_A$ meeting $B_A$ in $\alpha_A$,
there could be an embedded, essential planar surface in $M_A$ having
several boundary components with slope $\alpha_A$ on $B_A$ and it
could be extended to an embedded punctured-disk in $M$; but its
boundary slope on $B$ would be $\alpha$ and we would not get a new
slope.

So, suppose our embedded punctured-disk in $M$ does not meet
$\eta(A)$ in vertical annuli. Hence, there is just one component,
which by Lemma \ref{planar-pair-pants} is a punctured-annulus having
one of its boundary components in $B$ with slope $\gamma$.

We consider the two possibilities as to where the other boundary
component of this punctured-annulus is; it is either in $B''$ or in
$B_A$.

If the other boundary component is in $B''$, then the
punctured-annulus meets $B_A$ in a slope, say $\mu_A'$ that can be
Dehn twisted about $\alpha_A$ to $\mu_A$, the meridian on $B_A$.
Furthermore, $\mu_A'$ must bound a punctured-disk in $M_A$.  Thus
given any slope $\beta$ on $B$, we can construct a punctured-annulus
in $\eta(A)$ having one boundary on $B$ with slope $\beta$, the
other boundary on $B''$, and punctures in $\mu_A$. We can then Dehn
twist this punctured-annulus about an essential annulus in $\eta(A)$
between $B''$ and $B_A$ that has slope $\alpha''$ on $B''$ and slope
$\alpha_A$ on $B_A$. This gives a punctured-annulus in $\eta(A)$
meeting $B$ in slope $\beta$ and punctures in $B_A$ all having slope
$\mu_A'$. It follows that $\beta$ bounds a punctured-disk in $M$.
This gives the conclusion for part 2(b).

Finally, we assume our embedded punctured-disk in $M$ does not meet
$\eta(A)$ in a vertical annulus and there are no embedded
punctured-disks in $M_A$ meeting $B_A$ in a slope obtained by
twisting $\mu_A$ about $\alpha_A$.

Now, by Lemma \ref{planar-pair-pants}, we have that our
punctured-disk must meet $\eta(A)$ in a punctured-annulus having one
boundary the slope $\gamma$ in $B$, the other boundary in $B_A$ and
its punctures in $B''$. Hence, the slope on $B_A$ can be any
obtained from $\gamma_A$ by a Dehn twist about $\alpha_A$.

Clearly, if we have a punctured-disk in $M_A$ with boundary in
$B_A$, we can use Lemma \ref{planar-pair-pants} to construct a
punctured-disk in $M$ with boundary in $B$ and the slopes will
satisfy the relationship of our conclusions in our lemma.\end{proof}

\begin{remark} From Part 2(c) of Theorem \ref{planar-induct-step} we
see that if there is an essential annulus $A$ between $B$ and $B''$,
then for every slope $\gamma_A$ in $B_A$ that bounds a
punctured-disks in $M$, with punctures in $\bdy M_A\setminus B_A$,
there is generated an infinite family of punctured-disks in $M$ each
with its  boundary in $B$ having slopes corresponding to Dehn
twisting about $A$.\end{remark}

\vspace{.125 in}\noindent-{\it Algorithm for Condition $E_\gamma$.}
Our proof in the case of Condition $E_\gamma$ is by induction on the
number $n$ of boundary components of the link-manifold $M$. We note
that in the proof of Theorem \ref{find-planar} in the previous
section, we did not use induction but we did use Dehn fillings to
regularly reduce the number of boundary components. We could not use
induction in the previous section because we may, after filling, no
longer have an essential planar surface.  Now, for an embedded
punctured-disk, as in Conditions $E_\gamma$ and $E$, the existence
of an embedded punctured-disk assure the existence of an essential
one.

\begin{lem} \label{punctured-essential} Suppose $M$ is an irreducible
link-manifold, $B$ a component of $\bdy M$ and $\gamma$ a slope in
$B$. If $M$ contains an embedded punctured-disk with boundary having
slope $\gamma$ and punctures in $M\setminus B$, then $M$ contains an
embedded, essential punctured-disk with boundary having slope
$\gamma$ and punctures in $\bdy M\setminus B$.\end{lem}
\begin{proof} We have observed that a properly embedded surface in a
link-manifold is essential if and only if it is incompressible and
not an annulus or torus parallel into the boundary. Thus suppose $P$
is a properly embedded punctured-disk  with $bd(P)$ having slope
$\gamma$ in $B$ and all punctures in $\bdy M\setminus B$. If $P$ is
incompressible, then since it has precisely one boundary component
in $B$, it must be essential; however, if $P$ were not
incompressible, a compressions on $P$ leaves an embedded planar
surface with precisely one component of its boundary ($bd(P)$) in
$B$ and all other punctures in $\bdy M\setminus B$. After a finite
number of compressions, we get the desired, essential
punctured-disk.\end{proof}

\begin{thm}{\bf (Condition $E_\gamma$)} \label{cond-E-gamma} Given a
link-manifold $M$, a component $B$ of $\bdy M$, and  a slope
$\gamma$ in $B$, there is an algorithm to decide if $M$ contains an
embedded, essential punctured-disk with boundary having slope
$\gamma$ and punctures in $\bdy M\setminus B$. If there is one, the
algorithm will construct one.\end{thm}
\begin{proof} Our proof is by induction on the number $n$ of
boundary components of the given link-manifold. We begin with one
boundary component.

\vspace{.125 in}\noindent{\bf ($\boldsymbol{n=1}$)} $M$ is a
knot-manifold.

Suppose the knot manifold $M$ is given by the triangulation $\T$.
For a knot manifold Condition $E_\gamma$ is the question of whether
there is an embedded disk with boundary slope $\gamma$. Consider the
normal solution space $\S_\gamma(M,\T)$ of normal surfaces
satisfying the boundary condition of meeting $\bdy M$ in the slope
$\gamma$. Compute the fundamental surfaces for $M$ in
$\S_\gamma(M,\T)$.  From Lemma \ref{find-disk}, there is an
embedded, essential disk with boundary having slope $\gamma$ in $M$
if and only if there is one among these fundamental surfaces.

\vspace{.125 in}\noindent{\bf ($\boldsymbol{n \Rightarrow n+1}$)}
Our induction hypothesis is that given a link-manifold $M$, a
component $B$ of $\bdy M$, and a slope $\gamma$ on $B$, then if $M$
has no more that $n, n\ge 1$, boundary components, we can decide if
there is an embedded, essential punctured-disk in $M$ with boundary
the slope $\gamma$ and punctures in $\bdy M\setminus B$.
Furthermore, if there is one, the algorithm will construct one.

So suppose we are given a link-manifold $M$ with $n+1$ boundary
components; $B$ is a selected boundary component of $\bdy M$; and
$\gamma$ is a slope on $B$.

By Theorem \ref{prime-decomp}, we can construct a prime
decomposition of $M = p(S^2\times S^1)\#q(\rppp)\#r(D^2\times
S^1)\#M_1\#\cdots\#M_K$, where each $M_i$ is given by a
$0$--efficient triangulation. In the construction of such a prime
decomposition, if $B$ results in being a boundary component of some
solid torus in the decomposition, then the algorithm actually
constructs an embedded, essential disk with boundary in $B$. From
Proposition \ref{planar-prime-factor} and the fact that there is a
unique slope in the boundary of a solid torus bounding an embedded,
essential disk, we can decide Condition $E_\gamma$. Furthermore, if
$B$ is in a prime factor having fewer boundary components than $M$,
then we can invoke our induction hypothesis.

Hence, we may assume $M$ has $n+1$ boundary components and is given
by a $0$--efficient triangulation $\T$. It follows that M has no
normal $2$-spheres, has no embedded $\rpp$, the only normal disks
are vertex-linking; hence, $M$ is irreducible and
$\bdy$--irreducible. Furthermore, all vertices of $\T$ are in $\bdy
M$ and each component of $\bdy M$ has precisely one vertex. We
consider the solution space $\S_\gamma(M,\T)$ and compute its
fundamental surfaces.

Case A. Suppose there is a fundamental surface that is an essential
annulus $A$ having boundary in distinct components $B'\ne B''$ of
$\bdy M$.

If $B'=B\ne B''$, then $A$ is the desired punctured-disk having
boundary $\gamma$ and one puncture in $B''$.

If $B'\ne B\ne B''$, then split the link-manifold $M$ along $A$. We
get a new link-manifold $M_A$. By part 1(b) of Lemma
\ref{planar-induct-step} there is an embedded punctured-disk in $M$
having boundary slope $\gamma$ in $B$ if and only if there is one in
$M_A$. The link-manifold $M_A$ has $n$ boundary components. By our
induction hypothesis, we can decide if $M_A$ contains an embedded
punctured-disk with boundary having slope $\gamma$ in $B$. If we
find one in $M_A$, then Lemma \ref{planar-induct-step} tells us how
to construct one in $M$.

We remark that $M_A$ is not given by a triangulation but a
cell-decomposition. However, from this cell-decomposition, we can
construct a link-manifold $N\subset M_A$ (but not necessarily $M_A$)
given by a triangulation, typically with fewer tetrahedra than $\T$;
$N$ has no more than $n$ boundary components, one of which is $B$,
and has the property that there is an embedded punctured-disk in $N$
with boundary  slope $\gamma$ in $B$ if and only if there is an
embedded punctured-disk with boundary slope in $M_A$ with boundary
slope $\gamma$.

Finally, we are left with the situation that $M$ is given by a
$0$--efficient triangulation $\T$ and there are no embedded,
essential annuli having their boundary in distinct components of
$\bdy M$ and are also represented in the solution space
$\S_\gamma(M,\T)$.

If there is any embedded punctured-disk in $M$ with boundary having
slope $\gamma$, say $P$, then, having the hypothesis of Proposition
\ref{ALE-link-bound}, there exists a constant $C$, depending only on
$M$ and $\T$, so that the average length of the punctures in $P$ is
no larger than $3C$, $\lambda^{-\gamma}_{av}\le 3C$.

Construct the set $\{\alpha_1,\ldots,\alpha_K\}$ of all slopes on
the components of $\bdy M\setminus B$ that have the property
$L(\alpha_i)\le C$. Now, construct the triangulated Dehn fillings,
$M(\alpha_1), \ldots, M(\alpha_K)$. There is a punctured-disk in $M$
with boundary having slope $\gamma$ and punctures in $\bdy
M\setminus B$ if and only if there is an embedded punctured-disk in
some $M(\alpha_i)$ with boundary having slope $\gamma$ and punctures
in $\bdy M(\alpha_i)\setminus B$. But each $M(\alpha_i)$ is a
link-manifold with no more than $n$ boundary components. By our
induction hypothesis, we can decide if there is such an embedded
punctured-disk in some $M(\alpha_i)$. If there is one, the algorithm
will construct one; the component of such a punctured-disk that
meets $B$ after removing the interior of the layered solid torus
$\bbb{T}(\alpha_i)$ is a desired solution for the given manifold
$M$.\end{proof}

\vspace{.125 in}\noindent-{\it Algorithm for Condition $E$.} Recall
that if we designate a slope $\mu$ on the boundary of a
link-manifold, we often call it a meridian and any slope $\lambda$
with  $\langle\mu,\lambda\rangle =1$ is called a longitude (with
respect to $\mu$). Condition $E$ ask if given a meridian, can we
decide if a longitude bounds a punctured-disk. {\it A priori} this
seems much more daunting that answering Condition $E_\gamma$, as we
are asking here if some slope among an infinite family of slopes
bounds a punctured-disk. We are able to solve this problem quite
easily and, in fact, use  Condition $E_\gamma$ to answer Condition
$E$. We first need an observation.

Suppose $M$ is a link-manifold, $B$ is a component of $\bdy M$, and
$\mu$  is a slope in $B$. Let $\eta(B) = B\times [0,1]$ be a small
product neighborhood of the boundary component $B$ and choose
notation so $B = B\times 1$. Let $\mu_0 = \mu\times \{0\}$ (here we
pick any representative for $\mu$). If we remove the interior of a
small tubular neighborhood of $\mu_0$, we have a link-manifold with
one more boundary component than $M$. We  use $M^\mu$ to denote this
new link-manifold and use $B^\mu$ to designate the boundary
component of $M^\mu$ that is not a component of $\bdy M$. The unique
slope on $B^\mu$ that bounds a disk in the tubular neighborhood of
$\mu_0$ is called the meridional slope on $B^\mu$ and designated by
$\mu^*$; it is a meridian in the typical sense as it bounds the
meridional disk of the tube about $\mu_0$. The embedded, essential
annulus $A = \mu\times [0,1]\cap M^\mu$  has boundary slope $\mu$ in
$B$ and boundary slope a longitude, say $\lambda^*$ (with respect to
$\mu^*$), in $B^\mu$. We say $M^\mu$ is the link-manifold {\it
obtained from $M$ by Dehn drilling along the slope $\mu$}.

\begin{prop}\label{dehn-drill} Let $M$ be a link-manifold, $B$ a component of $\bdy
M$,
and $\mu$ a slope in $B$. There is an embedded punctured-disk in $M$
with boundary the slope of a longitude in $B$ if and only if every
longitude in $B$ bounds a punctured-disk in $M^\mu$.\end{prop}

\begin{figure}
\psfrag{m}{$\mu$}\psfrag{n}{$\mu^*$}\psfrag{o}{$\lambda^*$}
\psfrag{l}{$\lambda$}\psfrag{b}{$B^\mu$}\psfrag{L}{$\lambda+\mu$}\psfrag{B}{\large{$B$}}
\psfrag{P}{\large{$P$}}
\psfrag{Q}{\large{$P^\mu$}}\psfrag{A}{\large{$A$}}{\epsfxsize = 3 in
\centerline{\epsfbox{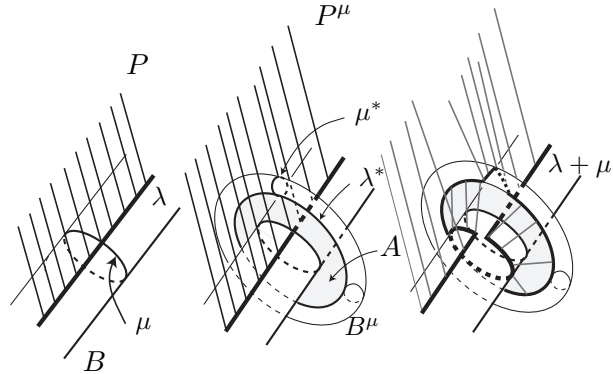}} } \caption{Dehn twisting a
punctured-disk with boundary a longitude in $B$ about the annulus
$A$ in the link-manifold $M^\mu$ obtained by Dehn drilling along the
slope $\mu$.} \label{f-dehn-drill}
\end{figure}

\begin{proof} Suppose there is a punctured-disk $P$  in $M$  having boundary a longitude
in $B$. If $P^\mu = P\cap M^\mu$, then (possibly after an isotopy)
we have $P^\mu$ meeting $B^\mu$ in the meridional slope $\mu^*$. See
Figure \ref{f-dehn-drill}. It follows that Dehn twisting $P^\mu$
about the annulus $A$ gives an embedded punctured-disk in $M^\mu$
with boundary a longitude in $B$; however, every longitude in $B$
can be realized this way.

Conversely, if there is an embedded punctured-disk in $M^\mu$ with
boundary a longitude in $B$, then there is one for every longitude
and, in particular, after an appropriate Dehn twist about $A$, there
will be one that meets $B^\mu$ in a single meridian. Such a
punctured-disk gives the desired punctured-disk in $M$ by Dehn
filling $M^\mu$ along the meridian in $B_\mu$ to reclaim
$M$.\end{proof}

\begin{thm}{\bf (Condition $E$)} Given a link-manifold $M$, a component $B$ of
$\bdy M$, and a slope $\mu$ in $B$, there is an algorithm to decide
if $M$ contains an embedded, essential punctured-disk with boundary
having slope of a longitude in $B$ (with respect to $\mu$) and
punctures in $\bdy M\setminus B$. If there is one, the algorithm
will construct one.\end{thm}
\begin{proof} This follows quite straight forward from
Proposition \ref{dehn-drill}. We construct the manifold $M^\mu$
obtained from $M$ by Dehn drilling along $\mu$. Select any longitude
$\lambda$ in $B$ and run the algorithm Condition $E_\lambda$. From
 Proposition \ref{dehn-drill}, there is an embedded punctured-disk in $M$ having
slope a longitude (not necessarily $\lambda$) and punctures in
$M\setminus B$ if and only if there is an embedded punctured-disk in
$M^\mu$ with boundary having slope $\lambda$ and punctures in $\bdy
M^\mu \setminus B$.\end{proof}

\begin{remark} We note that we have not said anything about a triangulation
of $M^\mu$. It may be the case that the slope $\mu$ is quite long
relative to the given triangulation on $M$, which indicates that we
may need to add a number of tetrahedra; also, we must add boundary
and at least one more vertex.  There is, however, a nice way to
triangulate $M^\mu$ that addresses these suspected needs in a very
visible way.

Suppose $M$ is given by a triangulation $\T$. Using the notation
from above, we see that a small neighborhood of $B\cup A\cup B^\mu$
is a pair of pants crossed with a circle. We shall add a ``pinched"
pair of pants crossed with the circle to $M$ along $B$ in such a way
as to get the manifold $M^\mu$ and extend $\T$ to a triangulation
$\T^\mu$ of the Dehn drilling $M^\mu$. See Figure
\ref{f-dehn-drill-triang}.

\begin{figure}
\psfrag{b}{$b$}
\psfrag{c}{$c$}\psfrag{d}{$b'$}\psfrag{x}{$x$}\psfrag{u}{\footnotesize{$u$}}\psfrag{v}{\footnotesize{$v$}}\psfrag{y}{$y$}\psfrag{z}{$z$}\psfrag{l}{$b'\times
S^1$}\psfrag{2}{$B^\mu = c\times
S^1$}\psfrag{4}{\footnotesize{$u\times S^1$}}\psfrag{3}{$b\times
S^1$}\psfrag{1}{$b'\times S^1$}\psfrag{R}{\large{$\td{R}$}}
\psfrag{S}{\large{${R}$}} \psfrag{T}{\large{$B^+ = R\times
S^1$}}\psfrag{i}{\small{identity}}{\epsfxsize = 3 in
\centerline{\epsfbox{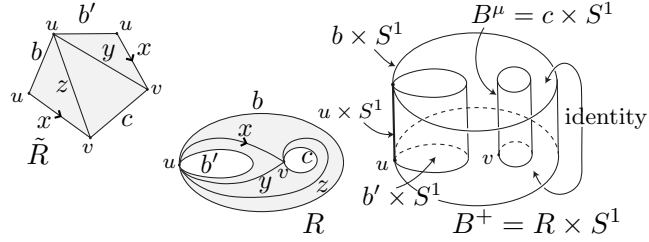}} } \caption{A
modification to the given triangulation for $M$ to a triangulation
of the manifold $M^\mu$ obtained by Dehn drilling along the slope
$\mu$ in $B$.} \label{f-dehn-drill-triang}
\end{figure}

Let $\td{R}$ be a pentagon in the plane with two boundary edges
labeled $x$ to be identified as indicated, three vertices labeled
$u$ to be identified, and two vertices labeled $v$ to be identified.
Label the three remaining boundary edges $b, b'$ and $c$, as shown.
Triangulate $R$ (any way is fine). After identification, this gives
a triangulation of a ``pinched" pair-of-pants, $R$. We take the
product $R\times S^1$, denoted $B^+$. We can triangulate $B^+$ with
nine tetrahedra so that $b,b',c, u\times S^1$ and $v\times S^1$ are
all edges.  We wish to identify the boundary $b'\times S^1$ of $B^+$
with $B$ to obtain the manifold $M^\mu$. To do this, the slope of
$u\times S^1$ must be identified with the slope $\mu$; then the
boundary $c\times S^1$ corresponds to $B^\mu$. We also want the
identification to be simplicial  so that we can extend the
triangulation $\T$ of $M$ to a triangulation of $M^\mu$.

To accomplish these requirements, we layer tetrahedra onto $B$, at
each layering getting a new triangulation of $M$, and eventually
arriving at a triangulation of $M$ that has added no new vertices
and has an edge, say $e^\mu$ with slope $\mu$. We identify $b'\times
S^1$ with $B$, taking $u\times S^1$ to $e^\mu$ and $b'$ to either of
the other edges. This determines a diagonal on the boundary
triangulation of $b'\times S^1$, which can be realized by a nine
tetrahedra triangulation of $B^+$.

This gives a minimal vertex-triangulation of $M^\mu$. Note that the
layerings on $B$ are determined by the length of $\mu$, a necessity
we expected, and the nine additional tetrahedra allow room for the
new boundary and vertex on $c\times S^1 = B^\mu$.
\end{remark}

\subsection{The family of slopes for punctured-disks.} From our methods
to obtain the algorithm for Conditions $E_\gamma$,  we are able, if
given a link-manifold and a component $B$ of $\bdy M$, to give a
construction that lists precisely those slopes on $B$ that are
boundary slopes of embedded punctured-disks in $M$.

It seems that algorithms for Condition $E_\gamma$ and Condition $E$
might follow from this result. However, the list may be infinite and
we have the classical problem in obtaining algorithms. Namely, if
there is an embedded punctured-disk in $M$ satisfying either of the
boundary conditions, $E_\gamma$ or $E$, then by systematically
constructing precisely those slopes in $B$ that can bound embedded,
punctured disks, we will encounter the corresponding slope and have
our solution. But if no slope does satisfy Condition $E_\gamma$ or
Condition $E$ and the set of possible slopes is infinite, when do we
stop looking? We believe, however, that the construction is quite
informative.

\begin{thm}\label{find-planar-bdd} Suppose we are given a link-manifold
$M$ and a component $B$ of $\bdy M$. There is a constructible set of
slopes $\S$ in $B$ so that
\begin{enumerate}
\item [(i)] any embedded punctured-disk in $M$ with boundary in $B$ and punctures in $\bdy M\setminus
B$ has boundary slope in $\S$; and
\item[(ii)] every slope in $\S$ is the boundary slope of an embedded punctured-disk in $M$ with
boundary in $B$ and punctures in $\bdy M\setminus B$.
\end{enumerate}
\end{thm}

\begin{proof} The proof of Theorem \ref{find-planar-bdd} is by induction on the
number, $n$, of boundary components of the given link-manifold $M$.

{\bf (n = 1)} For a knot-manifold, an embedded punctured-disk with
boundary in $B$ and punctures in $\bdy M\setminus B$ is an embedded
disk. By Lemma \ref{find-disk} there is an algorithm to decide if
$M$ contains an embedded disk with its boundary in $B$. If it does,
then $M$ splits as a connected sum of a solid torus with boundary
$B$ and a closed $3$--manifold. In this case $\S$ contains a single
slope. Otherwise, $\S=\emptyset$.

$\mathbf {(n\Rightarrow n+1)}$ Our induction step assumes that given
a link-manifold $M'$ that has no more than $n$ boundary components,
$n\geq 1$, along with a distinguished component $B'$ of $\bdy M'$,
then there is a constructible set of slopes $\S'$ in $B'$ so that
items $(i)$ and $(ii)$ are satisfied if we substitute $M', B'$ and
$\S'$ for  $M, B$ and $\S$, respectively.

Now, suppose we are given the link-manifold $M$ that has  $n + 1$
boundary components and $B$ is a specified component of $\bdy M$.

We use Theorem \ref{prime-decomp} to construct a prime decomposition
of $M$ so that $B$ is either the boundary of a solid torus or $B$ is
the boundary of a prime factor $M'$ given by a $0$--efficient
triangulation. In the former case $\S$ consists of the unique single
meridional slope; in the latter case, if $M'$ has fewer than $n+1$
boundary components, then we apply induction. Hence, we may assume
$M$ has $n+1$ boundary components and is given by a $0$--efficient
triangulation. Compute the fundamental surfaces of $M$.

We may assume no fundamental surface is an embedded, essential disk,
giving the following possibilities.

First, we consider the possibility that a fundamental surface is an
embedded, essential annulus having boundary in distinct components
of $\bdy M$.

We note that for a fixed component of $\bdy M$ there is at most one
slope that is in the boundary of an embedded, essential annulus; or
$M = S^1\times S^1\times I$ ($M$ has more than one boundary
component).

Suppose there is an embedded, essential annulus $A$ between the
components $B'$ and $B''$, $B'\ne B''$ of $\bdy M$.

There are two possibilities; either $B'\ne B\ne B''$ or $B'= B\ne
B''$.

{\bf 1. $B'\ne B\ne B''$.} Split $M$ at $A$ to get the link-manifold
$M_A$. Then $M_A$ has fewer boundary components than $M$ and $B$ is
a component of $\bdy M_A$. By Propositions \ref{planar-induct-step}
Part 1(b) we have that $\S = \S_A$, where $\S_A$ is the set of
slopes on $B$, as a component of $\bdy M_A$, that bound
punctured-disks in $M_A$ with punctures in $\bdy M_A\setminus B$.

To construct the slopes in $\S$,  we need to construct the slopes in
$\S_A$. In this situation, we provided a method for a triangulation
of $M_A$ in the Remark following the proof of Theorem
\ref{cond-E-gamma}.

{\bf 2. $B' = B (B\ne B``$).} Again, we split $M$ at $A$ to get the
link-manifold $M_A$. Then $M_A$ has fewer boundary components than
$M$ and $B_A$ becomes the distinguished component of $\bdy M_A$.

By Propositions \ref{planar-induct-step} we can have any one of
three possibilities. From Part 2(b) there are embedded
punctured-disks in $M$ with boundary meeting $B$ in every slope if
and only if there is an embedded punctured-disk in $M_A$ with
boundary meeting $B_A$ in a slope obtained from $\mu_A$ by Dehn
twisting about $\alpha_A$. In this case $\S$ consists of every slope
in $B$. From Part 2(c), if the previous situation does not hold,
then possibly the only embedded punctured-disk in $M$ meeting $B$ is
$A$ and there are no embedded punctured-disks in $M_A$ meeting
$B_A$. In this case $\S$ consists of precisely the slope $\alpha$.
Finally, we have the possibility of an embedded punctured-disk in
$M$ with boundary meeting $B$ in slope $\gamma$ if and only if there
is an embedded punctured-disk in $M_A$ with boundary meeting $B_A$
in a slope obtained from $\gamma_A$ by Dehn twisting about
$\alpha_A$. In this case, the slopes in $\S$ are generated by Dehn
twists about $\alpha_A$ of slopes in $\S_A$.

Again, to find the slopes in $\S$ we need to have the slopes in
$\S_A$, which begins with our having the manifold $M_A$ given to us
by a triangulation. As above we find this step described in the
Remark following the proof of Theorem \ref{cond-E-gamma}.

The only remaining situation is that no fundamental surface in
$(M,\T)$ is an essential disk or an essential annulus with its
boundaries in distinct components of $\bdy M$.

In this situation, if there is an embedded, essential
punctured-disk, say $P$, in $M$ with its boundary in $B$ and its
punctures in $\bdy M\setminus B$, then
$$P = \sum n_k F_k + \sum m_j K_j +\sum l_i A_i,$$ where $\chi(F_k)
< 0, K_j$ is a torus or Klein bottle, and $A_i$ is either an
annulus, with its boundary components in the same component of $\bdy
M$, or a M\"obius band. Thus we have a situation for the use of an
average length estimate. In this case, we can apply Proposition
\ref{ALE-link-bound} where we do not restrict the slope $\gamma$ and
we have  $b_\gamma = 1$. The average length estimate gives us that
there is a constant $C$ depending only on $M$ and the triangulation
$\T$ so that any planar surface we are interested in must have a
boundary shorter than $3C$.

Hence, there must be a short puncture (length less than $3C$) for
all embedded punctured-disk in $M$ having boundary $B$ and punctures
in $\bdy M\setminus B$.

Let $\{\alpha_1,\ldots, \alpha_N\}$ denote the slopes in the
components of $\bdy M\setminus B$ such that $L(\alpha_i)\le C$. If
there is any punctured-disk embedded in $M$ with boundary in $B$ and
punctures in $\bdy M\setminus B$, at least the set of punctures on
some boundary must have slope one of the $\alpha_i, 1\le i \le N$.

Let $(M(\alpha_i),\T(\alpha_i)),\ldots,(M(\alpha_N),\T(\alpha_N))$
be the triangulated Dehn fillings of $(M,\T)$ along the slopes
$\alpha_1,\ldots,\alpha_N$. Each $M(\alpha_j)$ is a link-manifold
having $n$ boundary components and $B$ is a component of $\bdy
M(\alpha_j)$ for all $j$. By our induction hypothesis, and for each
$j, 1\le j\le N$, there is a constructible set of slopes, which we
denote $\S(\alpha_j)$, in $B$ that satisfy conditions $(i)$ and
$(ii)$ of the theorem. In this case $\S
=\cup\S(\alpha_j)$.\end{proof}

Even though there are no embedded, essential annuli in $M$, an
$M(\alpha_j)$ can have an embedded, essential annulus between
distinct boundary components and possibly one with a component of
its boundary in $B$. If this were the case, then $\S(\alpha_j)$
could have infinitely many slopes along a line, leading to $\S$
having infinitely many slopes along a line. But since $M$, itself,
has no embedded, essential annulus with boundary in distinct
components of $\bdy M$ such a line of slopes in $\bdy M$ is {\it
not} obtained by Dehn twisting about an essential annulus in $M$.

While every slope in $\S$ is the boundary of an embedded
punctured-disk in $M$, there may be quite different punctured-disks
with the same slope for their boundary. Our methods do not find all
the possible punctured-disks.

\section{The word problem}

As mentioned earlier, this work was initiated by a study of singular
normal surfaces in an attempt to solve the Word Problem for
fundamental groups of $3$--manifolds. In the Introduction, we gave
the Word Problem for $3$--manifolds as a decision problem for
knot-manifolds; we then gave the Word Problem for finitely presented
groups as an analogous decision problem for link-manifolds.

We begin this section with the classical version of the Word Problem
for $3$--manifold groups and show that it is equivalent to the
version given in the Introduction.

\vspace{.125 in}\noindent {\footnotesize WORD PROBLEM} (closed
$3$--manifolds). {\it Given a closed $3$--manifold $M$ and a loop
$L$ in $M$. Decide if $L$ is contractible in $M$}.

It is quite straight forward that this version is equivalent to the
version that is in the Introduction; however, there is an
interesting point about the version from  the Introduction. Namely,
we again have the issue of having an infinity of longitudes from
which we want to know if one bounds a suitable punctured-disk. We
address this for singular punctured-disks in Lemma
\ref{sing-longitudes} below.

We repeat here for convenience the version given in the
Introduction.

\vspace{.125 in}\noindent {\footnotesize WORD PROBLEM} (closed
$3$--manifolds). {\it Given a knot-manifold $M$ and a meridian on
$\bdy M$. Decide if a longitude bounds a (possibly) singular
punctured-disk in $M$ with punctures meridians}.

We first show (briefly) why the two versions are equivalent. Let
$\eta(L)$ denote a small regular neighborhood of $L$ in $M$. Let
$M_L=M\setminus \open{\eta}(L)$; then $M_L$ is a knot-manifold and
we designate the unique slope on $\bdy M_L$ that bounds a disk in
$\eta(L)$ as the meridian slope, $\mu$ (it is a meridian in the
classical sense).

If $L$ bounds a singular disk, say $D$, in $M$, we can set $P=D\cap
M_L$ and we have a singular punctured-disk in $M_L$ with boundary a
longitude on $\bdy M_L$ and punctures meridians. Conversely, if a
longitude in $\bdy M_L$ bounds a singular punctured-disk in $M_L$
and punctures meridians, then $L$ bounds a singular disk in $M$.
Hence, the answer to one version of the Word Problem for
$3$--manifolds is yes if and only if the answer to the other version
is yes.

The issue about the infinity of longitudinal slopes is easier to
resolve in the singular case than in our earlier considerations
regarding embedded punctured-disks.

\begin{lem}\label{sing-longitudes} Given a knot-manifold $M$ and a meridian slope on its
boundary. If any longitude on $\bdy M$ bounds a singular
punctured-disk in $M$ with punctures meridians, then all longitudes
on $\bdy M$ bound a singular punctured-disk in $M$ with punctures
meridians.\end{lem}
\begin{proof} Consider Figure \ref{f-dehn-drill}. Let $\mu$
denote the meridional slope on $\bdy M$. Now, using the notation
from Figure \ref{f-dehn-drill}, we have that a longitude on $\bdy M$
bounds a singular punctured-disk with punctures meridians in $M$ if
and only if every longitude on $\bdy M$, considered as a component
of $\bdy M^\mu$ (the manifold obtained from $M$ by Dehn drilling
along $\mu$), bounds a singular punctured-disk in $M^\mu$ with
punctures meridians on $\bdy M$ or on $B^\mu$ (the other boundary
component of $M^\mu$).

So, suppose some longitude on $\bdy M$ bounds a singular
punctured-disk in $M$ with punctures meridians. Let $\gamma$ be any
longitude on $\bdy M$. Then $\gamma$ bounds a singular
punctured-disk, say $P^\mu_\gamma$, in $M^\mu$ with punctures
meridians on $\bdy M$ or on $B^\mu$. In fact, in this case there
must be precisely one puncture with boundary on $B^\mu$ and its
slope is a Dehn twist of the meridian $\mu^*$ about the longitude
$\lambda^*$ on $B^\mu$; and so has slope $\mu^*+k\lambda^*$ on
$B^\mu$ for some integer $k$. Again, see Figure \ref{f-dehn-drill}.
However, the slope $\mu^*+k\lambda^*$ on $B^\mu$  bounds a singular
punctured-disk in the solid torus $\eta(\mu_0)$ having $\abs{k}$
punctures  in $B^\mu$ with slope $\lambda^*$. The latter are all
isotopic to $\mu$ on $\bdy M$. Thus there is a singular
punctured-disk $D_k$ in $M$ with boundary $\mu^*+k\lambda^*$ on
$B^\mu$ and $\abs{k}$ punctures in $\bdy M$, each having slope
$\mu$. The punctured disk $P_k = P^\mu_k\cup D_k$ is a singular
punctured-disk with boundary $\gamma$ and punctures meridians in
$\bdy M$.

To change boundary slope, we need to change the number of punctures.
But if there is a singular punctured-disk in $M$ with boundary a
longitude and punctures meridians, then there is such a {\it
singular} punctured-disk in $M$ for {\it every} longitudinal
slope.\end{proof}

The classical statement of the Word Problem for finitely presented
groups typically has the following form.

\vspace{.125 in}\noindent {\footnotesize WORD PROBLEM} (finitely
presented groups). {\it Suppose the group $G$ is given by the finite
presentation $G=\langle X:R\rangle$. Given a word $w$ in the symbols
of $X$ decide if $w = 1$ in $G$}.

In the Introduction we gave the following decision problem for
link-manifolds as equivalent to the Word Problem for finitely
presented groups.

\vspace{.125 in}\noindent {\footnotesize WORD PROBLEM} (finitely
presented groups). {\it Given a link-manifold $M$, a component $B$
of $\bdy M$, and a meridian slope on $B$. Decide if there is a
(possibly) singular punctured-disk in $M$ with boundary slope a
longitude in $B$ and punctures  in $\bdy M\setminus B$ or meridians
on $B$}.

We have the following construction. Let $G=\langle X:R \rangle$ be a
finite presentation of the group $G$, where $X =\{x_1,\ldots,x_n\}$
and $R=\{r_1,\ldots,r_K\}$,  $r_j$ a word in the symbols of $X$. Let
$H_n = n(S^2\times S^1)=(S^2\times S^1)\#\cdots\#(S^2\times S^1)$,
where the right hand side has $n$ terms. $H_n$ is the closed
handlebody of rank $n$; its fundamental group is free of rank $n$
and we label a set of free generators by $x_1,\ldots,x_n$.

The relations $r_1,\ldots,r_K$ of $R$ can be represented by a
collection of pairwise disjoint embedded loops in $H_n$; we shall
denote these loops also by $r_1,\ldots,r_K$ in $H_n$. Let
$\eta(r_j), j=1,\ldots,K$, be pairwise disjoint small tubes about
$r_1,\ldots,r_K$, respectively. Let $M_G = H_n\setminus
\bigcup_{j=1}^K\open{\eta}(r_j)$ and denote the boundary components
of $M_G$ by $B_j = \bdy\eta(r_j)$.

Now, let $v_1,\ldots, v_K$ be $K$ distinct points and let $C_j =
v_j*B_j$ be the cone on $B_j$ with cone point $v_j$. Finally, let
$\hat{M}_G = M_G\bigcup_{j=1}^K C_j$.

The fundamental group of $\hat{M}_G$ is the group $G$.

Notice we can triangulate $\hat{M}_G$ getting an ideal triangulation
of the interior of the compact link-manifold $M_G$,
$\hat{M}_G\setminus \{v_1,\ldots,v_K\} = \open{M}_G$; we call the
vertices $v_1,\ldots,v_K$ ideal vertices.

Now, suppose $w$ is a word in the symbols of $X$. Then $w$ can be
represented by an embedded loop, also denoted $w$, in $\hat{M}_G$
missing the vertices $v_1,\ldots,v_K$. Hence, $w$ can be represented
by a loop in $M_G$ that is equivalent to $w$ in $\hat{M}_G$. We
shall continue to call this loop $w$ as a loop in $M_G$. Let
$\eta(w)$ be a small regular neighborhood of $w$ in $M_G$ and set $M
= M_G\setminus \open{\eta}(w)$. Let $B = \bdy \eta(w)$. Then $M$ is
a link-manifold with a distinguished boundary component $B$ and a
natural meridional slope, say $\mu$.

\begin{lem}The word $w=1$ in $G$ if and only if the loop $w$ is
contractible in $\hat{M}_G$ if and only if for the meridional slope
$\mu$ on $B$ there is a singular punctured-disk in $M$ with boundary
a longitude in $B$ and punctures meridional in $B$ or on $\bdy
M\setminus B$.\end{lem}

We conclude with a curious set of relations. If we write Condition
$E$ as $E(1)$, in the case of one boundary component, then $E(1)$ is
the classical Knot-Problem; namely, is a given knot the unknot. If
we write $S(1)$ for the Word Problem for $3$-manifolds, then $S(1)$
is the singular version of the embedded version $E(1)$. So, the Word
Problem for $3$--manifolds is the singular version of the Knot
Problem.

Now,  write $E(n)$ for the link-manifold version of Condition $E$
with $n$ boundary components and call it the Link Problem.
Similarly, write $S(n)$ for the $3$--manifold version of the Word
Problem for finitely presented groups having $n-1$ relations, $n>1$.
Then the Word Problem for finitely presented groups is the Link
Problem for $3$--manifolds. A decision problem for $3$--manifolds
without a general solution.

We note that in using  embedded normal surface theory, we employed
induction to go from $E(1)$ to the general solution $E(n)$; however,
following the announced solution of the Geometrization Conjecture,
we have $S(1)$ solvable; so, there is no way to go from $S(1)$ to
$S(n)$, leaving us with believing that singular normal surface
theory can not be built as a direct analogue to the embedded theory.
Of course, there are many other reasons to draw this conclusion.

\end{document}